\documentclass[10pt]{amsart}
\usepackage{amsmath}
\usepackage{amssymb}
\usepackage{amsfonts}
\usepackage{amsthm}
\usepackage{enumerate}
\usepackage[all]{xy}
\usepackage[latin1]{inputenc}
\usepackage{graphicx}
\usepackage{latexsym}
\usepackage{accents}
\usepackage{mathdots}
\usepackage{mathrsfs}
\usepackage{color}
\input xy
\usepackage{tikz}
\usetikzlibrary{arrows}
\usetikzlibrary{decorations,arrows}
\usetikzlibrary{decorations.markings}
\usetikzlibrary{arrows.meta, bending}
\usetikzlibrary{positioning}

\usepackage[notref,notcite]{showkeys}

\makeatletter
\newtheoremstyle{definition}
        {5pt}
        {3pt}
        {}
        {0pt}
        {\scshape}
        {.}
        {5pt}
        {\thmname{#1} \thmnumber{#2} \thmnote{[#3]}} 

\newtheoremstyle{theorems}
        {5pt}
        {3pt}
        {\itshape}
        {0pt}
        {\scshape}
        {.}
        {5pt}
        {\thmname{#1} \thmnumber{#2}\thmnote{[#3]}} 

\swapnumbers \theoremstyle{theorems}

\newtheorem{Theo}{Theorem}[subsection]
\newtheorem{Prop}[Theo]{Proposition}
\newtheorem{Cor}[Theo]{Corollary}
\newtheorem{Lemma}[Theo]{Lemma}

\newtheorem{Theo(Aus)}[Theo]{Theorem (Auslander)}
\newtheorem{Prop(BG)}[Theo]{Proposition (Bongartz-Gabriel)}
\newtheorem{Lemma(Asashiba)}[Theo]{Lemma(Asashiba)}
\newtheorem{Lemma(Gab)}[Theo]{Lemma(Gabriel)}
\newtheorem{Theo(Mil)}[Theo]{Theorem (Milicic)}

\theoremstyle{definition}
\newtheorem{Defn}[Theo]{Definition}

\newtheorem{Defn(Asashiba)}[Theo]{Definition (Asashiba)}
\newtheorem{Remark}[Theo]{Remark}
\newtheorem{Exam}[Theo]{Example}

\theoremstyle{theorems}

\newcommand{\cA}{\mathcal{A}}

\newcommand{\cP}{\mathcal{P}}
\newcommand{\cI}{\mathcal{I}}

\def\Sa{\hbox{${\mathit\Sigma}$}}
\def\Oa{\hbox{${\mathit\Omega}$}}
\def\Ta{\hbox{${\mathit\Theta}$}}
\def\La{\hbox{$\it\Lambda$}}

\newcommand{\Hom}{{\rm Hom}}

\newcommand{\Ext}{{\rm Ext}}
\newcommand{\End}{{\rm End}}

\newcommand{\mmod}{{\rm mod}}
\newcommand{\soc}{{\rm soc \hspace{.4pt}}}
\newcommand{\rad}{{\rm rad \hspace{.4pt}}}

\def\Mod{\hbox{{\rm Mod}{\hskip 0.3pt}}}

\def\ModLa{\hbox{{\rm Mod}{\hskip 0.5pt}\La}}

\def\mf{\mathfrak}

\newcommand{\mk}{\mathfrak}

\newcommand{\m}{\hspace{-.5pt}}

\begin{document}

\title[Almost split sequences]{\sc Auslander-Reiten theory via Nakayama duality in abelian categories \vspace{-5pt}}

\author[Z. Lin]{Zetao Lin \vspace{-6pt}} \author[S. Liu]{Shiping Liu \vspace{-6pt}
}





\address{Zetao Lin \\ Department of Mathematical Sciences, Tsinghua University, 100084, Beijing, P. R. China.}
\email{zetao.lin@foxmail.com}

\address{Shiping Liu \\ D\'epartement de math\'ematiques, Universit\'e de Sherbrooke, Sherbrooke, Qu\'ebec, Canada.}
\email{shiping.liu@usherbrooke.ca}

\maketitle

\centerline{\it Dedicated to the memory of Idun Reiten \vspace{-5pt}}

\begin{abstract}

Using the Nakayama duality induced by a Nakayama functor, we provide a novel and concise account of the existence of Auslander-Reiten dualities and almost split sequences in abelian categories with enough projective objects or enough injective objects. As an example, 
we establish the existence of almost split sequences ending with finitely presented modules and those starting with finitely copresented modules in the category of all modules over a small endo-local Hom-reflexive category. Specializing to algebras given by (not necessarily finite) quivers with relations, we further investigate when the categories of finitely presented modules, finitely copresented modules and finite dimensional modules have almost split sequences on either or both sides.

\vspace{5pt}

\end{abstract}

\section*{Introduction}

\vspace{2pt}

The Auslander-Reiten theory of irreducible maps and almost split sequences, introduced by Auslander and Reiten in the 1970's; see \cite{AuR2, MAIR}, is a powerful tool for studying categories of various kinds: most prominently, categories of modules over algebras and their derived categories; see \cite{ARS, GaR, Ha2, IRMVB}, and also, categories of Cohen-Macaulay modules over commutative algebras; see \cite{LW, RS}, categories of sheaves on varieties and schemes; see \cite{Aus3, Jor}, derived categories of differential graded algebras arising from topological spaces; see \cite{Jo1}, and most recently, 
extriangulated categories; see \cite{OIHY, LiN}. 

In this paper, we are mainly concerned with the existence of almost split sequences in an $R$-category $\mk A$, where $R$ is a commutative ring.
Write $D$ for the dual functor $\Hom_R(-, \hspace{-1pt} I_{\hspace{-.5pt}R}),$ where $I_{\hspace{-.5pt}R}$ is the minimal injective cogenerator for $R$-modules. Typically, \vspace{.5pt} one derives an almost split sequence 
$\xymatrixcolsep{18pt}\xymatrix{\hspace{-2pt} 0\ar[r] & \tau M\ar[r] & E \ar[r] & M\ar[r] & 0\hspace{-2pt}}$ from an Auslander-Reiten duality ${\rm Ext}_{\mk A}^1(X, \tau M) \cong D\underline{\mk A}(M, X)$ or ${\rm Ext}_{\mk A}^1(M, X) \cong D\overline{\mk A}(X, \tau M)$ for all $X\in \mk A$, or from their weak forms as described in \cite[(2.3)]{SLPC}, where $\tau$ is the AR-translation; see \cite{A2, ARS, GaR, LeZ}. The classical AR-translation for module categories is the dual of the transpose; see \cite[(I.3)]{A2}, \cite[(7.4)]{AIR}, \cite[(2.2)]{AuR2} and \cite[(1.6.1)]{RMV}. 
For lack of a transpose functor, we shall take a different 
approach. 

Let $\mathcal P$ be a strictly full additive subcategory of projective objects in $\mk A$. As defined in \cite[(5.4)]{LiN}, a Nakayama functor $\nu: \mathcal{P}\to \mk A$ induces a 
{\it Nakayama duality} $D\mk{A}(P, X)\cong \mk{A}(X, \nu P)$ for $P\in \mathcal P$ and $X\in \mk A$. 
It has been shown that a Nakayama duality induces Aualander-Reiten dualities in the derived category of $\mk A$; see \cite[(5.7)]{LiN}. Inspired by our recent work on the Auslander-Reiten dualities for graded mo\-dules; see \cite[(3.4)]{LL25}, we show that it also induces, in an elegant way, Auslander-Reiten dualities in $\mk A$. Indeed, defining 
the {\it AR-translate} $\tau_{_\delta}\hspace{-.8pt}M$ of an object $M$ with a projective presentation $\delta$ 
over $\mathcal P$; see (\ref{pseu_AR_trasl}), we show straightforwardly that $ {\rm Ext}_{\hspace{.8pt}\mk{A}\hspace{-.8pt}}^1(X, \tau_{_\delta\hspace{-.8pt}}M) \cong D\hspace{0.5pt}\underline{\mk{A}}(M,X)$ for $X\in \mk A$ \vspace{.8pt} if $\mk A$ has enough projective objects, and $D\hspace{0.05em}{\rm Ext}^1_{\hspace{.8pt}\mk{A}\hspace{-.8pt}}(M, X)\cong \overline{\mk{A}}(X, \tau_{_\delta\hspace{-.8pt}}M)$ for $X\in \mk A$ if $\mk A$ has enough injective objects; see (\ref{Naka-fcon}) and (\ref{AR-formu}). Observe that the second isomorphism induces only a monomorphism ${\rm Ext}^1_{\hspace{.8pt}\mk{A}\hspace{-.8pt}}(M, X)\to D\overline{\mk{A}}(X, \tau_{_\delta\hspace{-.8pt}}M),$ not an isomorphism as stated above. 

In order to derive almost split sequences from Auslander-Reiten dualities or their weak forms, one requires a key property that the AR-translation preserves strong indecomposability. For this purpose, we assume that $\mathcal P$ is Hom-reflexive and Krull-Schmidt; see (\ref{KS_section}). In this case, we obtain a {\it Nakayama equivalence} $\nu: \mathcal{P}\to \nu \mathcal{P};$ see \cite[(5.6)]{LiN} and an AR-translation uniquely defined (up to isomorphism) for objects finitely presented over $\mathcal P$; see (\ref{AR-trsl}). These lead directly to the desired property; see (\ref{tau1}), without involving a transpose functor as the classical approach; see \cite[(7.4)]{AIR} and \cite[(IV.1)]{ARS}. This enables us to apply the result in \cite[(2.3)]{SLPC} to establish the existence of almost split sequences ending with objects finitely presented over $\mathcal P$ if $\mk A$ has enough projective objects, and almost split sequences starting with objects finitely copresented over $\nu \mathcal{P}$ if $\mk A$ has enough injective objects; see (\ref{C_mod_ass}).

Our results will be applicable to many concrete abelian categories, including potentially categories of sheaves on topological spaces. In Section 3, for example, we shall consider the category $\Mod \mathscr{C}$ of all left modules over a small endo-local Hom-reflexive $R$-category $\mathscr{C}$. By the Yoneda Lemma, the full subcategory ${\rm proj}\mathscr{C}$ of finitely generated projective mo\-dules in $\Mod \mathscr{C}$ is Hom-reflexive and Krull-Schmidt. By constructing a Nakayama functor $\nu: {\rm proj}\mathscr{C}\to \Mod\mathscr{C};$ see (\ref{NK_functor}), we obtain an existence theorem for almost split sequences in $\Mod\mathscr{C}$ ending with finitely presented modules and for those starting with finitely copresented modules; see (\ref{AR-seq}). Note that the first part of this existence theorem was stated without proof or reference by Auslander in \cite{A2}.
However, the most relevant result that we can find in the existent literature is the existence of the Auslander-Reiten duality, and consequently, the existence of almost split sequences, in the category of finitely presented $\mathscr{C}$-modules when $\mathscr{C}$ is a dualizing $R$-variety; see \cite[(7.4)]{AIR}. It is not clear how this approach works for $\Mod \mathscr{C},$ whether $\mathscr{C}$ is a dualizing $R$-variety or not. 

In section 4, we shall specialize to an algebra $\La$ given by a quiver with relations. In case $\La$ is locally semiperfect; see (\ref{lspa}), which can be viewed as a small endo-local Hom-reflexive category, we obtain immediately a general existence theorem for almost split sequences in $\Mod\La$; see (\ref{ass_lfd}). So we will focus on the subcategories $\mmod^+\hspace{-3pt}\La,$ $\mmod^-\hspace{-3pt}\La$ and $\mmod^{\hspace{.5pt}b \hspace{-3pt}}\La$ of finitely presented modules, finitely copresented modules and finite dimensional modules in $\Mod\La$, respectively. In case $\La$ is locally semiprimary given by locally finite quivers; see (\ref{lspri_def}), the almost split sequences in $\mmod^+\hspace{-3pt}\La$ and $\mmod^-\hspace{-2.5pt}\La$ are almost split in $\Mod\La$; see (\ref{ass_fpres}) and (\ref{ass_fcp}). And we find conditions for these subcategories to have almost split sequences on either or both sides; see (\ref{ass_fp_cat}), (\ref{ass_fcp_cat}) and (\ref{ass_fdm_cat}). In particular, we obtain examples in which $\mmod^{\hspace{.5pt}b \hspace{-3pt}}\La$ has almost split sequences, but it has neither enough projective objects nor enough injective objects; see (\ref{ex_ass_fdim_cat}). In the hereditary case, our results strengthen the corresponding ones in \cite{BLP}. We conclude 
by emphasizing the broad significance of the representation theory of algebras defined by infinite quivers with or without relations: it allows us to apply Bongartz and Gabriel's covering theory to study representations of finite dimensional algebras; see \cite{BauL, BaL, BoG, DSk}; and it establishes deep connections to other domains such as coalgebras and comodules; see \cite{CKQ, CDi}, non-commutative algebraic geometry; see \cite{IRMVB} and cluster theory; see \cite{HJor,LiP}.

\section{Preliminaries}

The objective of this section is to fix some terminology and notation, which will be used throughout this paper, and collect some preliminary results. 

\subsection{\sc Reflexive modules} Throughout this paper, $R$ stands for a commutative ring with an identity. We denote by $\Mod R$ the category of $R$-modules and by $I\!_R$ a minimal injective cogenerator for $\Mod R$. A significant role in our investigation will be played by the {\it dual functor} $D=\Hom_R(-, I\!_R): \Mod R\to \Mod R$ will play. Given an $R$-module $M,$ there exists a canonical $R$-linear monomorphism $\sigma_{\hspace{-1.5pt}_M}: M\to D^2\hspace{-1pt}M$, sending $m$ to the evaluating  at $m$ function. We shall say that $M$ is {\it reflexive} if $\sigma_{\hspace{-1.5pt}_M}$ is bijective. The following statement is known; see \cite[(1.3)]{LiN}.

\begin{Lemma}\label{Reflex}

Let $R$ be a commutative ring. Then the category $\mmod R$ of reflexive $R$-modules is abelian, and the dual functor restricts to a duality $D\hspace{-1.5pt}:\hspace{-1.5pt} \mmod R \hspace{-1.5pt}\to \hspace{-1.5pt} \mmod R.$
    
\end{Lemma}

\begin{Remark}\label{sp_rmod} (1) The category $\mmod R$ contains all $R$-modules of finite length.

\noindent (2) In case $R$ is complete noetherian local, $\mmod R$ contains all noetherian $R$-modules and all artinian $R$-modules.

\end{Remark}

\vspace{-5pt}

\subsection{\sc Krull-Schmidt categories}\label{KS_section} An $R$-category is a category in which the morphism sets are $R$-modules and the composition of morphisms is $R$-bilinear. We shall compose morphisms in any $R$-category from right to left. A full subcategory of an $R$-category is called {\it strictly full} if it is closed under isomorphisms. All functors between $R$-categories are assumed to be $R$-linear. 

Let $\cA$ be an $R$-category. Given objects $X, Y\in \cA$, we write $\cA(X, Y)$ for the $R$-module of morphisms $f: X\to Y.$ A non-zero object $X$ in $\cA$ is called {\it strongly indecomposable} if $\cA(X, X)$ is local; and {\it Krull-Schmidt} if it is a finite direct sum of strongly indecomposable objects. 
We shall say that $\cA$ is {\it Krull-Schmidt} if every non-zero object is Krull-Schmidt,  and {\it Hom-reflexive} (respectively, {\it Hom-noetherian, 
Hom-finite}) if the $R$-modules $\mathcal{A}(X, Y)$ are reflexive (respectively, noetherian, 
of finite length), for all objects $X, Y$ in $\cA$. Note that a Hom-finite $R$-category is Hom-reflexive and Hom-noetherian. 
Moreover, if $R$ is complete noetherian local, then a Hom-noetherian $R$-category is Hom-reflexive.

\begin{Prop}\label{Ku-sum}

Let $\cA$ be a strictly full additive subcategory of an abelian $R$-category $\mk A$. If $\cA$ is Krull-Schmidt, then it is closed under direct summands in $\mk A$. The converse holds if $\cA$ is either Hom-finite or Hom-noetherian with $R$ being complete local noetherian.
    
\end{Prop}

\noindent{\it Proof.} Firstly, assume that $\cA$ is Krull-Schmidt. Let $X\in \cA$ with a nonzero direct summand $Y$ in $\mk A$. Then, $X=X_1\oplus \cdots \oplus X_n$, where $X_i\in \cA$ with $\mathcal{A}(X_i,X_i)$ being local. Since $\mathfrak{A}(X_i, X_i)=\cA(X_i, X_i)$, we see that $Y$ is a direct sum of some of the $X_i$; see \cite[(1.1)]{SLPC}. In particular, $Y\in \cA$. 

Conversely, assume that $\cA$ is closed under direct summands in $\mk A$. Then,  the idempotent endomorphisms in 
$\cA$ split. Assume further that $\cA$ is either Hom-finite or Hom-noetherian with $R$ being complete noetherian local. Then, for any $X\in \cA$, the $R$-algebra $\mathcal{A}(X, X)$ is either artianian; or noetherian with $R$ being complete noetherian local. In either case, $\mathcal{A}(X, X)$ is semi-perfect; see \cite[(21.35)]{Lam}. Thus, $\cA$ is Krull-Schmidt; see \cite[(1.1)]{SLPC}. The proof of the proposition is completed.

\subsection{\sc Minimal morphisms} Let $\cA$ be an $R$-category. 
An epimorphism $f: X\to Y$ in $\cA$ is called {\it superfluous} if any morphism $g: M\to X$ in $\cA$ such that $fg$ is an epimorphism is an epimorphism. Dually, a monomorphism $f: X\to Y$ is called {\it essential} if any morphism $h: Y\to N$ such that $hf$ is a monomorphism is a monomorphism. More importantly, a morphism $f: X\to Y$ in $\cA$ is called {\it left minimal} if every morphism $g: Y\to Y$ such that $gf=f$ is an automorphism; and {\it right minimal} if every morphism $h: X\to X$ such that $fh=f$ is an automorphism. 
For later reference, we quote the following statement from \cite[(1.4)]{HKS}.

\begin{Prop}\label{min_decom}

Let $\cA$ be an additive $R$-category. 

\begin{enumerate}[$(1)$]

\vspace{-1.5pt}

\item Every nonzero morphism \vspace{.8pt} $f: X\to Y$ in $\cA,$ where $Y$ is Krull-Schmidt, can be decomposed as 
$f=(f_1, 0)^T: X\to Y_1 \oplus Y_2,$ where $f_1: X\to Y_1$ is left minimal. 

\vspace{1pt}

\item Every nonzero morphism $g: X\to Y$ in $\cA,$ where $X$ is Krull-Schmidt, can be decomposed as $f=(g_1, 0): X_1\oplus X_2 \to Y,$ where $g_1: X_1\to Y$ is right minimal. 

\end{enumerate} \end{Prop}

\vspace{-0pt}

\subsection{\sc Almost split sequence} Let $\cA$ be an exact $R$-category, that is, an extension-closed subcategory of an abelian $R$-category $\mk A$. An object $X$ in $\cA$ is called {\it ext-projective} \vspace{-1pt} if every short exact sequence $\hspace{-2pt} \xymatrixcolsep{18pt}\xymatrix{0\ar[r] & Z\ar[r] & Y \ar[r] & X\ar[r] & 0\hspace{-2pt} }$ in $\cA$ splits; and 
{\it ext-injective} if every short exact sequence $\hspace{-1pt} \xymatrixcolsep{18pt}\xymatrix{0\ar[r] & X\ar[r] & Y \ar[r] & Z\ar[r] & 0\hspace{-2pt}}$ in $\cA$ splits. If $\cA$ is abelian, then the projective objects in $\cA$ coincide with the Ext-projective objects, and the injective objects coincide with the Ext-injective objects. 

\vspace{1pt}

Let  $f: X\to Y$ be a morphism in $\cA$. One says that $f$ is {\it irreducible} if $f$ is neither a section nor a retraction, and given any factorization $f=gh$, either $h$ is a section or $g$ is a retraction. Moreover, $f$ is called {\it left almost split} if $f$ is not a section and every non-section morphism $g:X\to Z$ factors through $f$; and {\it minimal left almost split} if it is left minimal and left almost split. Dually, one defines {\it right almost split morphisms} and {\it minimal right almost split morphisms} in $\cA$. Furthermore, an {\it almost split sequence} in $\cA$ is a short exact sequence \vspace{-5pt} 
$$\xymatrixcolsep{18pt}\xymatrix{0\ar[r] & X\ar[r]^f & Y \ar[r]^g & Z\ar[r] & 0,} \vspace{-4pt}$$ where $f$ is minimal left almost split and $g$ is minimal right almost. In such an almost split sequence, $X$ is the {\it starting term} and $Z$ is the {\it ending term}. 

\vspace{1pt}

We shall say that $\mathcal{A}$ has almost split sequence {\it on the right} if every strongly indecomposable and not ext-projective object is the ending term of an almost split sequence, and $\mathcal{A}$ has almost split sequences  {\it on the left} if every strongly indecomposable and not ext-injective object is the starting term of an almost split sequence. Finally, we shall say that $\mathcal{A}$ has almost split sequences if it has almost split sequences on the right and on the left.

\subsection{\sc Stable categories} Let $\cA$ be an exact $R$-category. A morphism $f\hspace{-1.5pt}: \hspace{-1.5pt} X \hspace{-1.5pt}\to \hspace{-1.5pt} Y$ in $\cA$ is called {\it projectively trivial} if it factors through every proper epimorphism $g: M\to Y$ in $\cA,$ and {\it injectively trivial} if it factors through every proper monomorphism $h: X\to N$ in $\cA$. In case $\cA$ is abelian with enough projective objects, a morphism in $\cA$ is projectively trivial if and only if it factors through a projective object. And in case $\cA$ is abelian with enough injective objects, a morphism in $\cA$ is injectively trivial if and only if it factors through an injective object. 

\vspace{1pt}

Write $\mathscr{P}$ and $\mathscr{I}$ for the ideals in $\cA$ generated by the projectively trivial morphisms and the injectively trivial morphisms, respectively. The {\it projectively stable category} and the \vspace{.5pt} {\it injectively stable category} of $\cA$ are the quotient categories $\underline{\cA}=\cA/\mathscr{P}$ and $\overline{\hspace{-2.5pt}\cA\hspace{.5pt}}=\cA/\mathscr{I},$ respectively; see \cite[(9.2)]{GaR}, \cite[(2.1)]{LeZ}, \cite[Section 2]{SLPC}. 

\section{Auslander-Reiten Theory in abelian categories}

The main objective of this section is to provide a novel and concise account of the Auslander-Reiten theory in general abelian categories. First, we derive an Auslander-Reiten duality from the Nakayama duality induced by a Nakayama functor. In case the Nakayama functor is defined on a ``small" subcategory of projective objects, we show that the AR-translation is uniquely defined (up to isomorphism) on finitely presented objects and preserves the strong indeco\-mposability. This enables us to establish the existence of almost split sequences ending with finitely presented objects and those starting with finitely copresented objects.

\subsection{\sc Projective presentations and injective copresentations} Throughout this section, $\mk A$ stands for an abelian $R$-category. Let $\mathcal P$ and $\mathcal I$ be strictly full additive subcategories of projective objects and injective objects in $\mk A$, respectively. Consider an object $X$ in $\mk A$. A {\it projective cover} of $X$ {\it over $\mathcal{P}$} is a right minimal epimorphism $f: P\to X$ in $\mk A$ with $P\in \mathcal{P}$; and an {\it injective envelope} of $X$ {\it over $\mathcal{I}$} is a left minimal monomorphism $g: X\to I$ in $\mk A$ with $I\in \mathcal{I}$. Moreover, \vspace{-2.5pt} a {\it projective presentation} of $X$ {\it over $\mathcal P$} is an exact sequence $\hspace{-2pt}\xymatrixcolsep{20pt}\xymatrix{P_1\ar[r]^{d_1} & P_0 \ar[r]^{d_0}&X\ar[r]&0}$ in $\mk A$ with $P_1, P_0\in \mathcal P,$ which is {\it minimal} if $d_1, d_0$ are right minimal. \vspace{-2.5pt} Dually, a {\it injective copresentation} of $X$ {\it over $\mathcal I$} is an exact sequence 
$\hspace{-1.5pt}\xymatrixcolsep{20pt}\xymatrix{0\ar[r] & X\ar[r]^{d^{\hspace{.5pt}0}} & I^0 \ar[r]^{d^1} &I^1}\hspace{-3pt}$ in $\mk A$ with $I^0, I^1\in \mathcal I$, which is {\it minimal} if $d^{\hspace{.5pt}0}, d^1$ are left minimal. Clearly, projective cover, injective envelope, minimal projective presentation and minimal injective copresentation are unique up to isomorphism if they exist. We denote by $\mk{A}^+(\mathcal{P})$ and $\mk{A}^-(\mathcal{I})$ the full subcategories of objects in $\mk A$ with a projective presentation over $\mathcal{P}$ and with an injective copresentation over $\mathcal{I}$, respectively. 

\begin{Prop}\label{subcat_pres} Let $\mk A$ be an abelian $R$-category, and let $\mathcal P$ be a strictly full additive subcategory of projective objects in $\mk A$. Then, $\mk{A}^+(\mathcal{P})$ is closed under cokernels and extension-closed in $\mk A$. Moreover,

\begin{enumerate}[$(1)$]

\vspace{-1.5pt}

\item if $\cP\hspace{-.5pt}$ is Hom-reflexive $($respectively, Hom-noetherian, 
Hom-finite$\hspace{.8pt})$, then $\mk{A}^{\hspace{-.5pt}+\hspace{-1.5pt}}(\hspace{-.5pt}\mathcal{P}\hspace{-.3pt})$ is Hom-reflexive $($respectively, Hom-noetherian, 
Hom-finite$\hspace{.8pt});$

\vspace{.5pt}

\item if $\mathcal P$ is Krull-Schmidt, then  $\mk{A}^+\hspace{-1pt}(\hspace{-.5pt}\mathcal{P}\hspace{-.5pt})$ is closed under direct summands in $\mk A$, and every object in $\mk{A}^+\hspace{-1pt}(\hspace{-.5pt}\mathcal{P}\hspace{-.5pt})$ admits a minimal projective presentation over $\mathcal P.$

\end{enumerate} \end{Prop}

\noindent{\it Proof.} It follows from Proposition 2.1 in \cite{Aus1} that 
$\mk{A}^+(\mathcal{P})$ is closed under cokernels and extension-closed in $\mk A$.

(1) Let $X, X'$ be objects in $\mk{A}^+(\mathcal{P})$. \vspace{-3pt} By definition, there exist projective presentations $\hspace{-2pt}\xymatrixcolsep{20pt}\xymatrix{P_1\ar[r]^{d_1} &P_0 \ar[r]^{d_0}&X\ar[r]&0}$ and $\hspace{-2pt}\xymatrixcolsep{20pt}\xymatrix{P_1'\ar[r]^{d_1'} &P_0' \ar[r]^{d_0'}&X'\ar[r]&0}\hspace{-2pt}$ over $\cP$.
Set $\mathcal{M}(X, X'):=\{(f_1, f_0)\in \mk{A}(P_1, P_1') \oplus \mk{A}(P_1, P_1') \mid f_0d_1=d_1'f_1\}.$ \vspace{1.5pt} It is evident that we have an $R$-linear epimorphism 
$\varphi: \mathcal{M}(X, X')\to \mk{A}(X, X'),$ given by $(f_1, f_0)\mapsto f,$ where $f$ is the unique morphism making the diagram \vspace{-4pt}
$$\xymatrixrowsep{18pt}\xymatrix{P_1\ar[r]^{d_1} \ar[d]^{f_1}&P_0 \ar[r]^{d_0}\ar[d]^{f_0}&X\ar[d]^{f}\ar[r]&0 \\
P_1'\ar[r]^{d_1'} &P_0' \ar[r]^{d_0'}&X'\ar[r]&0}\vspace{-3pt}$$ commute. Suppose that $\cP$ is Hom-reflexive. Then $\mk{A}(P_1, P_1')$ and $\mk{A}(P_1, P_1')$ are reflexive $R$-modules. Since $\mmod R$ is abelian; see (\ref{Reflex}), $\mathcal{M}(X, X')$ is reflexive, and so is $\mk{A}(X, X').$ Similarly, if $\cP$ is Hom-noetherian 
or Hom-finite, then $\mk{A}^+(\mathcal{P})$ is Hom-noetherian 
or Hom-finite, respectively. 

\vspace{1pt}

(2) Assume that $\mathcal P$ is Krull-Schmidt. Then, \vspace{-4pt} $\mathcal P$ is closed under direct summands in $\mk A$. Let $X\in \mk{A}^+(\mathcal{P})$ with a projective presentation $\hspace{-3pt}\xymatrixcolsep{20pt}\xymatrix{P_1 \ar[r]^{d_1} &P_0\hspace{-1.5pt} \ar[r]^{d_{\hspace{.5pt}0}} & X\ar[r] &0}\hspace{-2pt}$ over $\mathcal P.$ By Proposition \ref{min_decom}, $d_{\hspace{.5pt}0}=(d'_{\hspace{.5pt}0}, 0):P_0=P'_0\oplus N_0\to X,$ where $d'_{\hspace{.5pt}0}$ is a right minimal epimorphism. Thus, $d'_0: P'_0\to X$ is an projective cover with $P'_0\in \mathcal{P}$. Moreover, \vspace{.5pt} ${\rm Im}(d_1)={\rm Ker}(d_{\hspace{.5pt}0})={\rm Ker}(d'_{\hspace{.5pt}0})\oplus N_0$. Thus, we have a decomposition $d_1=(f, 0)^T:P_1\to P'_0\oplus N_0,$ \vspace{.5pt} where $f: P_1\to P'_0$ is such that ${\rm Im}(f)={\rm Ker}(d'_0)$. By Proposition \ref{min_decom} again, \vspace{.5pt} $f=(d_1, 0): P'_1\oplus N_1\to P'_0,$ where $d'_1: P'_1\to P'_0$ is right minimal. It is easy to check that ${\rm Im}(d'_1)={\rm Im}(f)$. \vspace{-3pt} This yields a minimal projective presentation 
$\hspace{-4pt}\xymatrixcolsep{20pt}\xymatrix{P'_1 \ar[r]^{d'_1} &P'_0\hspace{-1.5pt} \ar[r]^{d'_{\hspace{.5pt}0}}& X\ar[r] & 0 \hspace{-2pt}}$ over $\mathcal P.$ 

\vspace{.5pt}

Suppose now that $X=X_1\oplus X_2$ with canonical projections $p_i: X\to X_i.$ In particular, we obtain epimorphisms $p_id'_0: P'_0\to X_i,$ for $i=1, 2$. \vspace{-2pt} As seen above, there exist short exact sequences 
$\hspace{-2pt}\xymatrixcolsep{18pt}\xymatrix{0\ar[r]& K_i \ar[r]^{f_i} & L_i \ar[r]^{g_i} & X_i \ar[r] & 0,}\hspace{-2pt}$ where $L_i\in \mathcal{P}$ and $g_i$ is a projective cover, for $i=1, 2$. By the uniqueness of projective cover, we see easily that $K_1\oplus K_2\cong {\rm Ker}(d'_0)={\rm Im}(d'_1).$ \vspace{-3pt} Therefore, we have epimorphisms $h_i: P_1'\to K_i$, for $i=1, 2.$ \vspace{-2pt} This yields projective presentations
$\xymatrixcolsep{20pt}\xymatrix{\hspace{-2pt} P'_1 \ar[r]^-{f_ih_i} & L_i \ar[r]^{g_i} & X_i \ar[r] & 0}\hspace{-2pt}$ over $\mathcal{P},$ for $i=1, 2$. So, $X_1, X_2\in \mk A^+(\mathcal P)$. The proof of the proposition is completed.

\vspace{3pt}

Dually, we have the following statement.

\begin{Prop}\label{subcat_copres} Let $\mk A$ be an abelian $R$-category, and let $\mathcal I$ be a strictly full additive subcategory of injective objects in $\mk A$. Then, $\mk{A}^+(\mathcal{I})$ is  closed under kernels and extension-closed in $\mk A$. Moreover,

\begin{enumerate}[$(1)$]

\vspace{-2pt}

\item if $\mathcal I$ is Hom-reflexive $($respectively, Hom-noetherian, 
Hom-finite$\hspace{.8pt})$, then $\mk{A}^{\hspace{-.5pt}-\hspace{-1.5pt}}(\hspace{-.5pt}\mathcal{I}\hspace{-.3pt})$ is Hom-reflexive $($respectively, Hom-noetherian, 
Hom-finite$\hspace{.8pt});$

\vspace{1pt}

\item if $\mathcal I$ is Krull-Schmidt, then $\mk{A}^-(\mathcal{I})$ is closed under direct summands in $\mk A,$ and every object in $\mk{A}^-(\mathcal{I})$ admits a minimal injective copresentation over $\mathcal I.$ 

\end{enumerate} \end{Prop}


\vspace{-2pt}

\subsection{\sc Nakayama functor} The key ingredient in our approach to the Auslander-Reiten duality is a Nakayama functor defined as follows; see \cite[(5.4)]{LiN}. 

\begin{Defn}\label{N_functor}

Let $\mk A$ be an abelian $R$-category, and let $\mathcal{P}$ be a strictly full additive subcategory of projective objects in $\mk A$. \vspace{.5pt} A functor $\nu:\mathcal{P}\to \mk A$ is called a {\it Nakayama functor} \vspace{.5pt} provided, for all $P\in \mathcal{P}$ and $X\in \mk A,$ that there exists an $R$-linear isomorphism  $\phi_{P, X}: \mk{A}(X,\nu P)\to D\mk{A}(P,X),$ which is natural in $P$ and $X$.

\end{Defn}

\begin{Remark} \label{NF_im} Since the dual functor $D$ is exact, the image $\mathcal{I}_\nu$ of a Nakayama functor $\nu:\mathcal{P}\to \mk A$ is a strictly full additive subcategory of injective objects in $\mk A$.
    
\end{Remark}


\subsection{\sc Auslander-Reiten duality} In order to have an Auslander-Reiten duality, we need to introduce the following notion, which depends on not only the objects in $\mk A^+(\mathcal P)$ but also their projective presentations over $\mathcal P$.

\begin{Defn}\label{pseu_AR_trasl} Let $\mk A$ be an abelian $R$-category equipped with a Nakayama functor $\nu:\mathcal{P}\to \mk A$. \vspace{-4pt} Given an object $M$ in $\mk A^+(\mathcal P)$ with a projective presentation $\xymatrixcolsep{20pt}\xymatrix{\hspace{-2pt} \delta: P_1\ar[r]^-{d_1} & P_0\ar[r]^{d_0} & M \ar[r] & 0}\hspace{-2pt}$ over $\mathcal P$, we put $\tau_{_\delta}M= {\rm Ker}(\nu d_1),$ called the {\it Auslander-Reiten translate} of $M$ {\it associated with $\delta$}.
    
\end{Defn} 

The following statement is essential for the Auslander-Reiten duality.

\begin{Prop}\label{Naka-fcon}

Let $\mk A$ be an abelian $R$-category \vspace{-2pt} equipped with a Nakayama functor $\nu:\mathcal{P}\to \mk A$. Consider \vspace{-5pt} a projective presentation $\xymatrixcolsep{20pt}\xymatrix{\hspace{-2pt} \delta: P_1\ar[r]^-{d_1} & P_0\ar[r]^-{d_0} & M \ar[r] & 0}\hspace{-2pt}$ over $\mathcal P$ and a short exact sequence $\xymatrixcolsep{20pt}\xymatrix{0\ar[r] &X\ar[r]^f &Y\ar[r]^g &Z\ar[r] &0}\hspace{-2pt}$ in $\mk A.$ \vspace{.8pt} Then, we have an $R$-linear exact sequence
\vspace{-4pt}
$$\xymatrixcolsep{25pt}\xymatrix{0\ar[r] & \mk{A}(Z, \tau_{_\delta\hspace{-.8pt}} M)\ar[r]^{g^*} &\mk{A}(Y, \tau_{_\delta\hspace{-.8pt}} M)\ar[r]^{f^*} &\mk{A}(X, \tau_{_\delta\hspace{-.8pt}} M) \ar[r]^-\eta &&} \vspace{-6.5pt}$$ 
$$\hspace{-18pt}\xymatrixcolsep{25pt}\xymatrix{ & D\mk{A}(M,Z) \ar[r]^-{Dg_*} &D\mk{A}(M,Y) \ar[r]^-{Df_*} & D\mk{A}(M,X) \ar[r] &0.} \vspace{-0pt} $$ 

\end{Prop}

\noindent{\it Proof.} \vspace{-1pt} Applying the Nakayama functor $\nu$ to $\delta$, we obtain an injective copresentation $\xymatrixcolsep{20pt}\xymatrix{\hspace{-2pt} 0\ar[r] & \tau_{_\delta\hspace{-.8pt}}M \ar[r] & \nu P_{1} \ar[r]^-{\nu d_{1}} &\nu P_{0}}$ over $\mathcal I_\nu$. Given an object $N$ in $\mk A$, applying $\mk{A}(N, -)$ to this injective copresentation yields an exact sequence \vspace{-5pt} $$(*) \qquad \xymatrixcolsep{25pt} \xymatrix{0\ar[r] & \mk{A}(N, \tau_{_\delta\hspace{-.8pt}}M) \ar[r] & \mk{A}(N, \nu P_{1}) \ar[r]^-{(\nu d_1)_*} & \mk{A}(N, \nu P_{0}).} \vspace{-3pt}$$ 

On the other hand, applying $D\mk{A}(-,\! N)$ to $\delta$ and using  the Nakayama duality, we obtain a commutative diagram \vspace{-4.5pt}
$$\xymatrixcolsep{28pt}\xymatrixrowsep{18pt}\xymatrix{
\mk{A}(N,\nu P_{1}) \ar[d]_{\cong} \ar[r]^-{(\nu d_1)_*} & \mk{A}(N,\nu P_0)\ar[d]_{\cong} \\
D\mk{A}(P_{1},N) \ar[r]^-{Dd_1^*} & D\mk{A}(P_{0},N) \ar[r]^-{Dd_0^*} 
& D\mk{A}(M,N) \ar[r] & 0,} \vspace{-1pt} $$ 
where the bottom row is exact. Combining this with $(*)$ yields an exact sequence \vspace{-3.5pt}
$$\xymatrixcolsep{25pt} \xymatrix{0\ar[r] & \mk{A}(N, \tau_{_\delta\hspace{-.8pt}}M) \ar[r] & \mk{A}(N, \nu P_{1}) \ar[r]^-{(\nu d_1)_*} & \mk{A}(N, \nu P_{0}) \ar[r] & D\mk{A}(M,N) \ar[r] & 0.}\vspace{-3pt}$$

Therefore, we obtain a commutative diagram with exact rows and columns
$$\xymatrixcolsep{25pt}\xymatrixrowsep{18pt}\xymatrix{&0 \ar[d] &0 \ar[d] &0 \ar[d]\\
0\ar[r] &\mk{A}(Z, \tau_{_\delta\hspace{-.8pt}}M) \ar[r]^{g^*} \ar[d] &\mk{A}(Y, \tau_{_\delta\hspace{-.8pt}}M) \ar[r]^{f^*} \ar[d] &\mk{A}(X,\tau_{_\delta\hspace{-.8pt}}M) \ar[d] \\
0\ar[r] &\mk{A}(Z,\nu P_{1}) \ar[r]\ar[d] &\mk{A}(Y,\nu P_{1}) \ar[r]\ar[d] &\mk{A}(X,\nu P_{1}) \ar[r]\ar[d] &0\\
0\ar[r] &\mk{A}(Z,\nu P_{0}) \ar[r]\ar[d] &\mk{A}(Y,\nu P_{0}) \ar[r]\ar[d] &\mk{A}(X,\nu P_{0}) \ar[r]\ar[d] &0\\
&D\mk{A}(M,Z) \ar[r]^{Dg_*} \ar[d] &D\mk{A}(M,Y) \ar[d]\ar[r]^{Df_*} &D\mk{A}(M,X) \ar[r]\ar[d]  &0, \\
&0 &0 &0
}$$ where the two middle rows are exact because $\nu P_{1}$ and $\nu P_0$ are injective. Applying the Snake Lemma, we obtain the desired exact sequence stated in the proposition. The proof of the proposition is completed.

\vspace{3pt}

We are ready to establish the Auslander-Reiten duality; compare \cite[(1.1), (1.2)]{HKra1}. Recall that $\underline{\mk{A}}=\mk{A}/\mathscr{P},$ the projectively stable category; and $\overline{\mk{A}}=\mk{A}/\mathscr{I},$ the injectively stable category, of $\mk A$.

\begin{Theo}\label{AR-formu}

Let $\mk A$ be an abelian $R$-category equipped with a Nakayama functor $\nu:\mathcal{P}\to \mk A$. \vspace{-2.6pt} Consider an object $M$ in $\mk{A}^+(\mathcal P)$ with a projective presentation $\xymatrixcolsep{20pt} \xymatrix{\hspace{-4pt} \delta: \; P_1\ar[r]^-{d_1} & P_0\ar[r]^{d_0} & M \ar[r] & 0} \hspace{-3pt}$ over $\mathcal P.$

\begin{enumerate}[$(1)$]

\vspace{-1pt}

\item If $\mk A$ has enough projective objects, \vspace{1pt} then there exists an $R$-linear isomorphism $D\hspace{0.5pt}\underline{\mk{A}}(M,X)\cong {\rm Ext}_{\hspace{.8pt}\mk{A}\hspace{-.8pt}}^1(X, \tau_{_\delta\hspace{-.8pt}}M)$ for  $X\in \mk A$, which is natural in $X$.

\vspace{0.5pt}

\item If $\mk A$ has enough injective objects,  \vspace{1pt}  then there exists an $R$-linear isomorphism $\overline{\mk{A}}(X, \tau_{_\delta\hspace{-.8pt}}M)\cong D\hspace{0.05em}{\rm Ext}^1_{\hspace{.8pt}\mk{A}\hspace{-.8pt}}(M, X)$ for  $X\in \mk A$, which is natural in $X$.

\end{enumerate} \end{Theo}

\noindent{\it Proof.} (1) Suppose that $\mk A$ has enough projective objects. \vspace{-1.5pt} Fix an object $X\in \mk A$. We have a short exact sequence 
$\xymatrixcolsep{20pt}\xymatrix{\hspace{-2pt} 0\ar[r] & L \ar[r]^-{q} &P\ar[r]^{p} & X \ar[r] &0,}\hspace{-1pt}$ where $P$ is a projective object in $\mk A$. It is easy to see that there exists an $R$-linear exact sequence \vspace{-10pt}
$$\xymatrix{ 0\ar[r] & D\underline{\mk{A}}(M,X) \ar[r] & 
D\mk{A}(M,X) \ar[r]^-{Dp_*} & D\mk{A}(M, P).} \vspace{-2pt} $$ 

Thus, we have an $R$-linear isomorphism $D\underline{\mk{A}}(M,X)\cong{\rm Ker} (D p_*)$, which is clearly natural in $X$. On the other hand, applying $\mk{A}(-, \tau_{_\delta\hspace{-.8pt}}M)$ to the above short exact sequence yields an $R$-linear exact sequence \vspace{-4pt} $$\xymatrixcolsep{18pt}\xymatrix{0\ar[r] &\mk{A}(X, \tau_{_\delta\hspace{-.8pt}}M) \ar[r]^-{p^{*}} &\mk{A}(P, \tau_{_\delta\hspace{-.8pt}}M) \ar[r]^-{q^{*}} & \mk{A}(L, \tau_{_\delta\hspace{-.8pt}}M)
\ar[r] &{\rm Ext}_{\hspace{.8pt}\mk{A}\hspace{-.8pt}}^1(X, \tau_{_\delta\hspace{-.8pt}}M) \ar[r] &0.} \vspace{-5pt}$$ 

Thus, we obtain an $R$-linear isomorphism ${\rm Coker}(q^{*})\cong {\rm Ext}_{\hspace{.8pt}\mk{A}\hspace{-.8pt}}^1(X, \tau_{_\delta\hspace{-.8pt}}M),$ which is also natural in $X$. Further, by Proposition \ref{Naka-fcon}, there exists an exact sequence \vspace{-6pt}
$$\xymatrixcolsep{25pt}\xymatrix{\mk{A}(P, \tau_{_\delta\hspace{-.8pt}} M)\ar[r]^{q^*} &\mk{A}(L, \tau_{_\delta\hspace{-.8pt}} M) \ar[r]^-\eta & D\mk{A}(M, X) \ar[r]^{Dp_*} & D\mk{A}(M,P)}\vspace{-5pt}$$ 

This yields an $R$-linear isomorphism 
${\rm Ker}(Dp_*) = {\rm Im}(\eta) \!\cong\! {\rm Coker}(q^{*}),$
which is natural in $X$. \vspace{1pt} Combining the above isomorphisms, we obtain an $R$-linear isomorphism 
$ D\hspace{0.05em}\underline{\mk{A}}(M,X)  \! \cong \!\Ext_{\hspace{.8pt}\mk{A}\hspace{-.8pt}}^1(X, \tau_{_\delta\hspace{-.8pt}}M),$ which is natural in $X$.

\vspace{1pt}

(2) \vspace{-1pt} Suppose that $\mk A$ has enough injective objects. \vspace{-1pt} Fix an object $X\in \mk A$. We have a short exact sequence $\xymatrixcolsep{20pt}\xymatrix{\hspace{-1pt} 0\ar[r] & X \ar[r]^-q & I \ar[r]^-p & L\ar[r]&0,\hspace{-.5pt}}$ where $I$ is an injective object in $\mk A$. Then, we have an exact sequence \vspace{-6pt} 
$$\xymatrixcolsep{25pt}\xymatrix{\mk{A}(I,\tau_{_\delta\hspace{-.8pt}}M) \ar[r]^-{q^*} & 
\mk{A}(X,\tau_{_\delta\hspace{-.8pt}}M) \ar[r] & \overline{\mk{A}}(X,\tau_{_\delta\hspace{-.8pt}}M) \ar[r] &0.} \vspace{-5pt} $$ 

So, we have an $R$-linear isomorphism $\overline{\mk{A}}(X,\tau_{_\delta\hspace{-.8pt}}M) \cong {\rm Coker}(q^*),$ which is natural in $X.$ On the other hand, applying $\mk{A}(M,-)$ yields an $R$-linear exact sequence \vspace{-2pt} $$\xymatrixcolsep{18pt}\xymatrix{0\ar[r] &\mk{A}(M,X) \ar[r]^-{q_{*}} &\mk{A}(M,I) \ar[r]^-{p_{*}} &\mk{A}(M, L)
\ar[r] &\Ext_{\hspace{.8pt}\mk{A}\hspace{-.8pt}}^1(M, X) \ar[r] & 0.}  \vspace{-2pt}$$ 

This gives rise to an $R$-linear isomorphism $D\Ext_{\hspace{.8pt}\mk{A}\hspace{-.8pt}}^1(M,X)\cong {\rm Ker}(Dp_*)$, which is natural in $X$. 
Further, by Proposition \ref{Naka-fcon}, we have an $R$-linear exact sequence \vspace{-5pt}
$$\xymatrixcolsep{28pt}\xymatrix{\mk{A}(I,\tau_{_\delta\hspace{-.8pt}}M) \ar[r]^-{q^*} &\mk{A}(X,\tau_{_\delta\hspace{-.8pt}}M)\ar[r]^-\eta &D\mk{A}(M, L)\ar[r]^-{Dp_*} &D\mk{A}(M,I).}\vspace{-2pt}$$

This gives rise to an $R$-linear isomorphism ${\rm Ker}(Dp_{*})={\rm Im}(\eta)\cong {\rm Coker}(q^*),$ which is natural in $X$. \vspace{1pt} Combining the above $R$-linear isomorphisms, we obtain an $R$-linear isomorphism 
${\rm \overline{\mk{A}}}(X, \tau_{_\delta\hspace{-.8pt}}M) 
\cong D{\rm Ext}_{\hspace{.8pt}\mk{A}\hspace{-.8pt}}^1(M, X),$ which is natural in $X$. The proof of the theorem is completed.


\subsection{\sc Auslander-Reiten translations} Let  $\nu: \mathcal{P} \to \mk A$ be a Nakayama functor. Given a strongly indecomposable object $M$ with a projective presentation $\delta$ over $\mathcal P$, in order to derive an almost split sequence from the Auslander-Reiten duality stated in Theorem \ref{AR-formu}(1), we need to ensure that $\tau_\delta M$ is strongly indecomposable. For this purpose, we assume that $\mathcal P$ is Hom-reflexive and Krull-Schmidt. In this setting, we have the following result; see 
\cite[(5.6)]{LiN}.

\begin{Prop}\label{fu-fai}

Let $\mk A$ be an abelian $R$-category equipped with a Nakayama functor $\nu:\mathcal{P}\to \mk A$, where $\mathcal P$ is Hom-reflexive and Krull-Schmidt. Then $\nu$ restricts to an equivalence $\nu: \mathcal{P} \to \mathcal{I}_\nu.$ In particular, $\mathcal{I}_\nu$ is Hom-reflexive and Krull-Schmidt.

\end{Prop}

\begin{Remark}

We call the equivalence $\nu: \mathcal{P}\to \mathcal{I}_\nu$ stated in Proposition \ref{fu-fai} the {\it Nakayama equivalence}, and we fix a quasi-inverse $\nu^-\!: \mathcal{I}_\nu\to \mathcal P$ for it.  
    
\end{Remark} 

By Propositions \ref{subcat_pres} and \ref{subcat_copres}, the exact categories $\mk{A}^+(\mathcal{P})$ and $\mk{A}^-(\mathcal{I}_\nu)$ are Hom-reflexive and closed under direct summands in $\mk A$. Moreover, every object in $\mk{A}^+(\mathcal{P})$ admits a minimal projective presentation over $\mathcal P,$ and every object in $\mk{A}^-(\mathcal{I}_\nu)$ admits a minimal injective copresentation over $\mathcal I_\nu$. This allows us to define, unique up to isomorphism, the {\it AR-translations} $\tau$ and $\tau^-$ as follows.

\begin{Defn}\label{AR-trsl}

Let $\mk A$ be an abelian $R$-category equipped with a Nakayama functor $\nu: \mathcal{P}\to \mk A$, where $\mathcal P$ is Hom-reflexive and Krull-Schmidt.

\begin{enumerate}[$(1)$]

\vspace{-6pt}

\item Given $M\in \mk{A}^+(\mathcal{P})$ of which \hspace{-4pt} $ \xymatrixcolsep{20pt}\xymatrix{P_1\ar[r]^{d_1} &P_0 \ar[r]^{d_0}&M\ar[r]&0}$ \hspace{-2pt} is a minimal projective presentation over $\mathcal{P}$, we define $\tau M={\rm Ker}(\nu d_1)$.

\vspace{-.5pt}

\item Given $N\in \mk{A}^-(\mathcal{I}_\nu)$ of which \hspace{-4pt} $\xymatrixcolsep{20pt}\xymatrix{0\ar[r] & N \ar[r]^{d^{\hspace{.5pt} 0}}&I_0\ar[r]^{d^1}&I_1}$ \hspace{-6pt} is a minimal injective copresentation over $\mathcal{I}_\nu$, we define $\tau^{-}\hspace{-1pt}N={\rm Coker} (\nu^-d^1)$. 

\end{enumerate} \end{Defn}


In order to show that $\tau$ and $\tau^-$ preserve the strong indecomposability, we slightly generalize a well-known fact stated in \cite[(II.4.3)]{A2}.

\vspace{-5pt}

\begin{Lemma}\label{LR-minimal} Let $\mk A$ be an abelian $R$-category, and let $\hspace{-3pt}\xymatrixcolsep{20pt}\xymatrix{X\ar[r]^{f} &Y \ar[r]^g &Z}\hspace{-3pt}$ be an exact sequence in $\mk A$.

\begin{enumerate}[$(1)$]

\vspace{-1.5pt}

\item If $\mk{A}(Z,Z)$ is local and $g$ is a nonsplit epimorphism, then $f$ is left minimal.

\vspace{1pt}

\item If $\mk{A}(X,X)$ is local and $f$ is a nonsplit monomorphism, then $g$ is right minimal.

\end{enumerate}\end{Lemma}

\noindent{\it Proof.} We shall only prove Statement (1). Suppose that $\mk{A}(Z,Z)$ is local and $g$ is a nonsplit epimorphism. \vspace{-2.8pt} Then, there exists a nonsplit short exact sequence
$\xymatrixcolsep{20pt}\xymatrix{\hspace{-2pt} 0\ar[r] & M\ar[r]^j \ar[r] & Y\ar[r]^{g} & Z\ar[r] & 0\hspace{-2pt}}$ in $\mk A$. By Lemma 4.3 in \cite[Chapter II]{A2}, $j$ is left minimal. And since $M={\rm Im}(f)$, there exists an epimorphism $v: X\to M$ such that $f=j v$. Assume that  $h f=f$ for some
$h: Y\to Y$. Then $hjv=jv,$ and hence, $hj=j$. Therefore, $h$ is an automorphism. 
The proof of the lemma is completed. 

\vspace{2pt}

The classical approach for showing that the Auslander-Reiten translation preserves indecomposability goes through the morphism category of projective modules and the stable categories; see \cite[(IV.1)]{ARS}. Here we provide a direct proof. 

\begin{Prop}\label{tau1} 

Let $\mk A$ be an abelian $R$-category equipped with a Nakayama functor $\nu: \mathcal{P}\to \mk A$, where $\mathcal P$ is Hom-reflexive and Krull-Schmidt.  Consider a strongly indecomposable object $M$ in $\mk A$.

\begin{enumerate}[$(1)$]

\vspace{-1.5pt}

\item If $M\in \mk{A}^+(\mathcal{P})$ is not projective, then $\tau M\in \mk{A}^-(\mathcal{I}_\nu)$ is strongly indecomposable and not injective such that $\tau^-(\tau M)\cong M$.

\vspace{0.5pt}

\item  If $M\in \mk{A}^-(\mathcal{I}_\nu)$ is not injective, then $\tau^-M\in \mk{A}^+(\mathcal{P})$ is strongly indecomposable and not projective such that $\tau(\tau^-M)\cong M$.

\end{enumerate}

\end{Prop}

\noindent{\it Proof.} We shall only prove Statement (1). \vspace{-3pt} Assume that $M\in \mk{A}^+(\mathcal{P})$ is nonprojective with a minimal projective presentation \hspace{-4pt} $\xymatrixcolsep{20pt}\xymatrix{P_1\ar[r]^{d_1} &P_0 \ar[r]^{d_0}&M\ar[r]&0}\hspace{-1pt}$ over $\mathcal{P}$. \vspace{-6pt} 
By definition,  $\tau M$ admits an injective copresentation \hspace{-4pt}
$\xymatrixcolsep{20pt}\xymatrix{0\ar[r] &\tau M\ar[r]^j &\nu P_1\ar[r]^{\nu d_1} &\nu P_0}\hspace{-3pt} $ over $\cI_\nu,$ which we claim is minimal. 
Indeed, $d_1$ is right minimal and $d_0$ is not a retraction. By Lemma \ref{LR-minimal}, $d_1$ is also left minimal. And by Proposition \ref{fu-fai}, $\nu d_1$ is left and right minimal. Since $\nu P_1$ is Krull-Schmidt, it follows from Proposition \ref{min_decom} that $j=(l, 0)^T: \tau M \to L\oplus U=\nu P_1,$ where $l: \tau M\to L$ is left minimal. Then, 
${\rm Ker}(\nu d_1)={\rm Im}(j)={\rm Im}(l).$ Hence, $\nu d_1=(f, g): \nu P_1=L\oplus U\to \nu P_0$, where $g: U\to \nu P_0$ is a monomorphism. Since $U$ is injective, $g$ is a section. And since $\nu: \mathcal{P}\to \mathcal{I}_\nu$ is an equivalence, we can write 
$d_1=(p,q): P_1=V\oplus W \to P_0$ such that $g=\nu(q: W\to P_0)$. Then, $q$ is a section such that $d_0q=0$. 
Since $d_0: P_0\to M$ is right minimal, we have $W=0$, and consequently, $U=0$. That is, $j$ is left minimal. Thus, $\tau M$ admits a minimal injective copresentation over $\mathcal{I}_\nu$ as claimed above. Applying $\nu^-$ to it, we see that $\tau^-(\tau M)\cong M$. 

Suppose that $\tau M$ is injective. Then, $j$ is a section. Being left minimal, $j$ is an isomorphism. Then, $\nu d_1=0$. Since $\nu d_1$ is left minimal, $\nu P_0=0$, and hence, $P_0=0,$ absurd. It remains to show that ${\mk A}(\tau M, \tau M)$ is local. Given any $f\in {\mk A}(\tau M, \tau M)$, we have a commutative diagram with exact rows \vspace{-5.5pt}  
$$\hspace{30pt}\xymatrixrowsep{20pt}\xymatrix{0\ar[r] &\tau M\ar[r]^j\ar[d]^{f} &\nu P_1\ar[r]^{\nu d_1}\ar[d]^{f_1} &\nu P_0\ar[d]^{f_0}\\0\ar[r] &\tau M \ar[r]^j &\nu P_1\ar[r]^{\nu d_1} &\nu P_0.}\vspace{-2pt} $$ Again since $\nu$ is an equivalence, we have a commutative diagram with exact rows \vspace{-4.5pt} 
$$\hspace{26pt}\xymatrixrowsep{20pt}\xymatrix{P_1\ar[r]^{d_1}\ar[d]^{g_1} &P_0\ar[d]^{g_0}\ar[r]^{d_0} &M \ar[r]\ar[d]^{g} &0 \\ P_1\ar[r]^{d_1} & P_0\ar[r]^{d_0} &M \ar[r] &0} \vspace{-3pt} $$ such that $\nu g_0=f_0$ and $\nu g_1=f_1$. Assume that $g$ is invertible. Since $d_0$ is superfluous; see \cite[(3.4)]{Krau} and $P_0$ is Krull-Schmidt, $g_0$ is an isomorphism. Similarly, $g_1$ is also an isomorphism. Thus, $f_0$ and $f_1$ are isomorphisms. So $f$ is invertible. Suppose that $g$ is not invertible. Since $\mk{A}(M, M)$ is local, $1\hspace{-1pt}_M-g$ is invertible. Observe that we obtain commutative diagrams from the above commutative diagrams by replacing $(f, f_1, f_0)$ and $(g_1, g_0, g)$ by $(1-f, 1-f_1, 1-f_0)$ and $(1-g_1, 1-g_0, 1-g)$, respectively. Using the same argument, we deduce that $1\hspace{-.4pt}_{\tau \hspace{-1pt} M}-f$ is invertible. So, $\mk{A}\hspace{-.4pt}(\tau\hspace{-1pt} M, \tau \hspace{-1pt} M)$ is local. The proof of the proposition is completed.

\vspace{2pt}

We conclude this subsection with some sufficient conditions for $\mk{A}^+(\mathcal{P})$ and $\mk{A}^-(\mathcal{I}_\nu)$ to be Krull-Schmidt.

\vspace{-2pt}

\begin{Lemma} \label{fpres_cat_KS}

Let $\mk A$ be an abelian $R$-category equipped with a Nakayama functor $\nu: \mathcal{P}\to \mk A$, where $\mathcal P$ is Krull-Schmidt. If $\mathcal P$ is either Hom-finite or Hom-noetherian with $R$ being complete noetherian local, then $\mk{A}^+(\mathcal{P})$ and $\mk{A}^-(\mathcal{I}_\nu)$  are Krull-Schmidt.

\end{Lemma}

\noindent {\it Proof.} Let $\mathcal{P}$ be Hom-finite. By Proposition \ref{fu-fai}, $\mathcal{I}_\nu$ is Hom-finite and Krull-Schmidt. 
By Propositions \ref{subcat_pres} and \ref{subcat_copres}, $\mk{A}^+(\mathcal{P})$ and $\mk{A}^-(\mathcal{I}_\nu)$ are Hom-finite and closed under direct summands in $\mk A$. By Proposition \ref{Ku-sum}, they are Krull-Schmidt. In case $R$ is complete noetherian local and $\mathcal{P}$ is Hom-noetherian, we can prove the result by the same argument. The proof of the lemma is completed.

\vspace{-2pt}

\subsection{\sc The existence theorem} We are ready to apply the result in \cite[(2.3)]{SLPC}
to derive almost split sequences from the Auslander-Reiten dualities in Theorem \ref{AR-formu}.

\vspace{-2pt}

\begin{Theo}\label{AR-seq}

Let $\mk A$ be an abelian $R$-category equipped with a Nakayama functor $\nu: \mathcal{P}\to \mk A$, where $\mathcal P$ is Hom-reflexive and Krull-Schmidt.

\begin{enumerate}[$(1)$]

\vspace{-1.5pt}

\item If $\mk A$ has enough projective objects, then it has an almost split sequence \vspace{-2.5pt}
$$\xymatrixcolsep{22pt}\xymatrix{0\ar[r] &\tau M\ar[r] & E \ar[r] & M \ar[r] &0,}\vspace{-5pt}$$ for every strongly indecomposable nonprojective object $M$ in $\mk{A}^+(\mathcal{P}).$

\vspace{1pt}

\item If $\mk A$ has enough injective objects, then it has an almost split sequence \vspace{-4.5pt}
$$\xymatrixcolsep{22pt}\xymatrix{ 0\ar[r] &N\ar[r] &E \ar[r] &\tau^{-}\hspace{-2.5pt} N \ar[r] &0,}\vspace{-6pt}$$  for every strongly indecomposable noninjective object $N$ in $\mk{A}^-(\cI_\nu).$ 

\end{enumerate} 

\end{Theo}

\noindent{\it Poof.} (1) Assume that $\mk A$ has enough projective objects. Consider a strongly indecomposable nonprojective $M\in\mk A^+(\cP)$. By Proposition \ref{tau1}(1), $\tau M$ is strongly indecomposable, and by Theorem \ref{AR-formu}(1), there exists a functorial isomophism $\Psi\hspace{-1pt}_M:\Ext^1_{\hspace{.4pt}\mk{A}\hspace{-.4pt}}(-,\tau M)\to D\hspace{0.05em}\underline{\mk{A}}(M,-).$ In particular, \vspace{.5pt} $\Ext^1_{\hspace{.4pt}\mk{A}\hspace{-.4pt}}(M,\tau M)\cong D\hspace{0.05em}\underline{\mk{A}}(M,M)$ as right $\mk{A}(M,M)$-modules. Since $\mk{A}(M,M)$ is local, 
we see that $D\hspace{0.05em}\underline{\mk{A}}(M,M)$ has a nonzero $\mk{A}(M,M)$-socle, and so does $\Ext^1_{\hspace{.4pt}\mk{A}\hspace{-.4pt}}(M,\tau M)$. By Theorem 2.3 in \cite{SLPC}, we obtain a desired almost split sequence as stated in Statement (1).

(2) Assume that $\mk A$ has enough injective objects. Consider a strongly indecomposable noninjective object $N\in \mk{A}^-(\cI_\nu)$. By Proposition \ref{tau1}(2), $\tau^-N$ is strongly indecomposable, and by Theorem \ref{AR-formu}(2), \vspace{.5pt} there exists a functorial isomorphism $\Phi_N: D^2\Ext^1_{\hspace{.4pt}\mk{A}\hspace{-.4pt}}(\tau^-\hspace{-1pt} N,-) \to D\hspace{0.05em}\overline{\mk{A}}(-, N)$. 
This induces a functorial monomorphism 
$\Psi_N: \Ext^1_{\hspace{.4pt}\mk{A}\hspace{-.4pt}}(\tau^{-}\!N,-)\to D\hspace{0.05em}\overline{\mk{A}}(-, N)$. Moreover, we have a left $\mk{A}(N,N)$-linear isomorphism $\Phi_{N, N}: D^2\Ext^1_{\hspace{.4pt}\mk{A}\hspace{-.4pt}}(\tau^{-}\!N, N) \to D\overline{\mk{A}}(N,N).$ \vspace{.5pt} By Propositions \ref{fu-fai} and \ref{subcat_copres}, $\mk{A}(N,N)$ is $R$-reflexive, \vspace{.5pt} and by Lemma \ref{Reflex}, both $D\overline{\mk{A}}(N,N)$ and $ D^2\Ext^1_{\hspace{.4pt}\mk{A}\hspace{-.4pt}}(\tau^{-}\!N, N)$ are $R$-reflexive. 
\vspace{.5pt} Therefore, $\Ext^1_{\hspace{.4pt}\mk{A}\hspace{-.4pt}}(\tau^{-}\!N, N) \cong D^2\Ext^1_{\hspace{.4pt}\mk{A}\hspace{-.4pt}}(\tau^{-}\!N, N) \cong D\overline{\mk{A}}(N,N)$ as left $\mk{A}(N,N)$-modules. Since $\mk{A}(N, N)$ is local, $D\overline{\mk{A}}(N, N)$ has a nonzero $\mk{A}(N, N)$-socle, and so does $\Ext^1_{\hspace{.4pt}\mk{A}\hspace{-.4pt}}(\tau^{-}\!N, N).$ By Theorem 2.3 in \cite{SLPC}, we obtain a desired almost split sequence as stated in Statement (2). The proof of the theorem is completed. 

\vspace{2pt}

As an easy example, we obtain Auslander's results stated in \cite[(II.6.3), (II.6.6)]{A2}.

\vspace{-2pt}

\begin{Theo(Aus)} Let $A$ be a noetherian $R$-algebra with an identity with $R$ complete noetherian local, and let $M$ be an indecomposable module in $\Mod A$.

\begin{enumerate}[$(1)$]

\vspace{-1.5pt}

\item If $M$ is \vspace{-.5pt} noetherian and not projective, then there exists an almost split sequence $\xymatrixcolsep{20pt}\xymatrix{\hspace{-2pt}0\ar[r] & \tau M \ar[r] & E \ar[r] & M\ar[r] & 0\hspace{-2pt}}$  in $\Mod A,$ where $\tau M$ is artinian.

\vspace{1pt}

\item If $M$ is \vspace{-1pt} artinian and not injective, then there exists an almost split sequence $\xymatrixcolsep{20pt}\xymatrix{\hspace{-2pt}0\ar[r] & M \ar[r] & E \ar[r] & \tau^-\hspace{-1pt} M\ar[r] & 0\hspace{-2pt}}$ in $\Mod A,$ where $\tau^-\hspace{-1pt}M$ is noetherian.

\vspace{-2pt}

\end{enumerate}\end{Theo(Aus)}

\noindent{\it Proof.} The category $\Mod A$ of left $A$-modules has enough projective and injective modules. And the full subcategories of noetherian modules and artinian mo\-dules in $\Mod A$ are Krull-Schmidt; see \cite[(I.5.1), (I.5.2)]{A2}. In particular,  
the full subcategory ${\rm proj}A$ of finitely generated projective modules in $\Mod A$ is Krull-Schmidt and Hom-noetherian; see \cite[(I.4.2)]{A2}. By Lemma 5.5 in \cite{LiN}, there exists a Nakayama functor $\nu=D\Hom_A(-, A): {\rm proj} A\to \Mod A.$ Further, a module in $\Mod A$ is noetherian if and only if it has a projective presentation over ${\rm proj} A;$ and artinian if and only if it has an injective copresentation over $\nu({\rm proj}A);$ see \cite[(I.5.2)]{A2}. Now, the statement follows from Theorem \ref{AR-seq}. The proof of the corollary is completed.

\vspace{-2pt}

\begin{Remark}

Let $\Sa$ be a ring with an identity. Auslander obtained an existence theorem for almost split sequences in $\Mod\Sa$ ending with finitely presented modules; see \cite[(II.5.1)]{A2}. To the best of our knowledge, there exists no such statement for almost split sequences in $\Mod\Sa$ starting with finitely copresented modules. 

\end{Remark}

\section{Applications to functor categories} 

An $R$-category is called {\it small} if its objects form a set and {\it endo-local} if all endomorphism algebras are local. The objective of this section is to apply our previous results to establish the existence of almost split sequences for finitely presented modules and for finitely copresented modules in the category of all modules over a small endo-local Hom-reflexive $R$-category. The key ingredients for this objective include a Nakayama functor and the sufficiency of projective and injective objects. The latter fact was stated in \cite{Wat}, but its proof seems to have never been published. 

\vspace{-4pt}

\subsection{\sc Modules over a category} Throughout this section, $\mathscr{C}$ stands for a small $R$-category. A {\it left $\mathscr{C}$-module} is a covariant functor $M: \mathscr{C}\to {\rm Mod} R.$ Given left $\mathscr{C}$-modules $M, N,$ a {\it $\mathscr{C}$-linear morphism} $f: M\to N$ consists of a family of $R$-linear maps $f_x: M(x)\to N(x)$ with $x\in \mathscr{C}$ such that $N(u)\circ f_x=f_y \circ M(u),$ for all morphisms $u: x\to y$ in $\mathscr{C}.$
Write $\Hom_{\mathscr{C}}(M, N)$ for the $R$-module of $\mathscr{C}$-linear morphisms $f: M\to N$. The category $\Mod\mathscr{C}$ of all left $\mathscr{C}$-modules is an abelian $R$-category with arbitrary products and coproducts; see \cite[Section 2]{A2}. A module $M\in \Mod\mathscr{C}$ is called {\it locally $R$-reflexive} if $M(x)$ is an $R$-reflexive for all $x\in \mathscr{C}.$ Since the category $\mmod R$ of reflexive $R$-modules is abelian, the full subcategory $\mmod \mathscr{C}$ of locally $R$-reflexive modules in $\Mod\mathscr{C}$ is abelian. 

We denote by $\mathscr{C}^\circ$ the opposite category of $\mathscr{C}:$ the objects are those of $\mathscr{C}$ and the morphisms $u^\circ: y\to x$ are induced by the morphisms $u: x\to y$ in $\mathscr{C}$. Given a module $M$ in $\Mod \mathscr{C}^\circ$, we define a module $\mf{D}M$ in $\Mod \mathscr{C}$ by $(\mf{D}M)(x)=D(M(x))$ and $(\mf{D}M)(u)=D(M(u^\circ))$ for all objects $x$ and morphisms $u$ in $\mathscr{C}.$ This yields an exact contravariant functor $\mf{D}: \Mod\mathscr{C}^\circ\to \Mod\mathscr{C}.$ 

\vspace{-2pt}

\begin{Prop}\label{LRMod}

Let $\mathscr{C}$ be a small $R$-category. The exact contravariant functor
$\mk{D}: \Mod\mathscr{C}^\circ \to \Mod\mathscr{C}$ restricts to a duality $\mk{D}: \mmod \mathscr{C}^\circ \to \mmod \mathscr{C}\hspace{-1.5pt}.$

\end{Prop}

\noindent{\it Proof.} Let $M\in \Mod\mathscr{C}^\circ.$ By Lemma \ref{Reflex}, there exist canonical $R$-linear isomorphisms $\sigma_{M, x}: M(x)\to D^2(M(x))$ with $x\in \mathscr{C}$, which form a $\mathscr{C}$-linear isomorphism, natural in $M$. The proof of the proposition is completed. 

\vspace{-5pt}

\subsection{\sc Projective modules} For each object $x$ in $\mathscr{C}$, we obtain a left $\mathscr C$-module $P_x=\mathscr{C}(x, -): \mathscr{C} \to \Mod R.$ And a morphism $u: y\to x$ in $\mathscr{C}$ induces a $\mathscr{C}$-linear morphism $\mathscr{C}(u, -): P_x\to P_y$ with $\mathscr{C}(u, -)_z = \mathscr{C}(u, z): \mathscr{C}(x, z)\to  \mathscr{C}(y, z)$ for $z\in \mathscr{C}$. We denote by $1_x$ the identity morphism associated with an object $x$ in $\mathscr{C}$.

\vspace{-1pt}

\begin{Prop}\label{Mor_proj_mod}

Let $\mathscr{C}$ be a small 
$R$-category. Given $M\in \Mod \mathscr{C}$ and $x\in \mathscr{C}$, we have an $R$-linear isomorphim, which is natural in $M,$ as follows$\,:$ \vspace{-4pt}
$$\Phi_{P_x, M}: \Hom_{\mathscr{C}}(P_x, M) \to M(x): f\mapsto f_x(1_x).\vspace{-4pt}$$

\end{Prop}

\noindent{\it Proof.} Fix $M\in \Mod \mathscr{C}$ and $x\in \mathscr{C}$. By the Yoneda Lemma, $\Phi_{P_x, M}$ is an $R$-linear isomorphism, which is natural in $M$ and $x$. If $\varphi: P_x\to P_y$ is a $\mathscr{C}$-linear morphism, by the Yoneda Lemma again, $\varphi=\mathscr{C}(u, -)$ for some $u: y\to x$ in $\mathscr{C}$. Thus $\Phi_{P_x, M}$ is natural in $P_x$ because it is natural in $x$.
The proof of the proposition is completed.


\begin{Remark}\label{proj_Cmod}

By Proposition \ref{Mor_proj_mod}, $P_x$ is projective in $\Mod\mathscr{C}\hspace{-1pt},$ for all $x\in \mathscr{C}$.

\end{Remark} 


\vspace{2pt}

The following result is stated with a sketched proof in \cite{Wat}. We include a short proof for the reader's convenience.

\begin{Prop}\label{Proj_C_mcat}

Let $\mathscr{C}$ be a small $R$-category. Then $\Mod\mathscr{C}$ has enough projective objects.

\end{Prop}

\noindent{\it Proof.} Fix $M\in \Mod\mathscr{C}$. Given $x\in \mathscr{C}$ with $m\in M(x),$ we write $x_m=x$. Then, \vspace{.5pt} $P=\textstyle\coprod_{x\in \mathscr{C}; \hspace{.5pt} m\in M(x)}\mathscr{C}(x_m, -)$ is projective in $\Mod \mathscr{C};$ see (\ref{proj_Cmod}). For $y\in \mathscr{C}$, \vspace{0pt} we have an $R$-linear map
$\psi_y: P(y)={\coprod}_{x\in \mathscr{C}; \hspace{.5pt} m\in M(x)}\mathscr{C}(x_m,y)\to M(y),$ defined by \vspace{.5pt} $\psi_y(u_{x,m})=M(u_{x,m})(m)$ for $u_{x,m}\in \mathscr{C}(x_m, y)$ such, for any $v:y\to z$ in $\mathscr{C}$, that \vspace{-8pt}
$$\xymatrixrowsep{20pt}\xymatrixcolsep{35pt}\xymatrix{\coprod_{{x\in \mathscr{C}; \hspace{.5pt} m\in M(x)}}\mathscr{C}(x_m,y) \ar@<7ex>[d]_{\coprod_{x\in \mathscr{C}; \hspace{.5pt} m\in M(x)}\mathscr{C}(x_m,v)} \ar[r]^-{\psi_y} & M(y) \ar[d]^{M(v)}\\ 
\coprod_{{x\in \mathscr{C}; \hspace{.5pt} m\in M(x)}}\mathscr{C}(x_m,z) \ar[r]^-{\psi_z} &M(z)}\vspace{-5pt}$$ commutes. Given $m\in M(y)$, we have $\hspace{.5pt}\psi_y(1_y)=M(1_y)(m)=m$. 
Thus, the $R$-linear epimorphisms $\psi_y$ with $y\in \mathscr{C}$ form a $\mathscr{C}$-linear epimorphism $\psi: P\to M.$ The proof of the proposition is completed.

\subsection{\sc Injective modules} Given $x\in \mathscr{C}$, we have a module $I_x=\mf{D}(P_x^\circ)$ in $\Mod\mathscr{C}\!,$ where $P_x^\circ=\mathscr{C}^\circ(x, -)=\mathscr{C}(-,x)\in \Mod \mathscr{C}^\circ.$ More explicitly, $I_x(y)=D\mathscr{C}(y,x)$ and  $I_x(u)=D\mathscr{C}(u,x): D\mathscr{C}(y, x) \to D\mathscr{C}(z,x),$ for objects $y$ and morphisms $u: y\to z$ in $\mathscr{C}.$ A morphism $v: x\to y$ in $\mathscr{C}$ induces a $\mathscr{C}$-linear morphism $D\mathscr{C}(-, v): I_y\to I_x$. And for $x\in \mathscr{C}$, we have an $R$-linear function $\theta_x: \mathscr{C}(x,x)\to I_R$ given by $f\mapsto f(1_x).$

\begin{Prop}\label{Mor-mod-inj}

Let $\mathscr{C}$ be a small $R$-category. Given  $M\in \Mod\mathscr{C}$ and  $x\in \mathscr{C},$ we have an $R$-linear isomorphism, which is natural in $M$, as follows$\,:$
\vspace{-1pt} $$\Psi_{M, I_x}: \Hom_{\mathscr{C}}(M, I_x)\to D(M(x)): f\mapsto \theta_x\circ f_x.\vspace{-0pt} $$ 
    
\end{Prop}

\noindent{\it Proof.} Fix $M\in \Mod\mathscr{C}$ and $x\in \mathscr{C}.$ Clearly, we have an $R$-linear map $\Psi_{M, I_x}$ as stated in the proposition. Assume that $f\in \Hom_{\mathscr{C}}(M, I_x)$ such that $\Psi_{M, I_x}(f)=0$. Given $y\in \mathscr C$ and $u\in \mathscr{C}(y,x),$ since $D\mathscr{C}(u,x) \circ f_y= f_x\circ M(u)$, we have \vspace{-2pt}
$$\begin{array}{rcl}
f_y(m)(u)&=& f_y(m)(\mathscr{C}(u,x)(1_x))\\
&=&\theta_x \left(D\mathscr{C}(u,x)(f_y(m))\right)\\
&=&\theta_x\left(f_x(M(u)(m))\right)\\
&=& \Psi_{M, I_x}(f)(M(u)(m))\\
&=& 0,   \vspace{-2pt}      
\end{array}$$ for all $m\in M(y)$ and $u\in \mathscr{C}(y,x).$ Thus, $f_y(m)=0$, for all $m\in M(y).$ That is, $f_y=0$ for all $y\in \mathscr{C}$, and hence, $f=0$. So, $\Psi_{M, I_x}$ is a monomorphism. Conversely, given any $R$-linear function $g: M(x)\to I_R$, we shall construct a $\mathscr{C}$-linear morphism $f: M\to I_x$. Indeed, for any $y\in \mathscr{C}$, we have an $R$-linear map \vspace{-2pt} $$\hspace{35pt} f_y: M(y) \to I_x(y)=D\mathscr{C}(y,x): m\mapsto f_y(m)\vspace{-1pt} $$ such that $f_y(m)(u)=g(M(u)(m))$, for $u \in  \mathscr{C}(y,x)$. Given $v \in \mathscr{C}(y,z),$ we claim that $D\mathscr{C}(v,x)\circ f_y=f_x \circ M(v)$.
Indeed, for $m\in M(y)$ and $u\in \mathscr{C}(z,x)$, we have \vspace{-1.5pt}
$$\hspace{-20pt} \begin{array}{rcl}
(D\mathscr{C}(v,x))(f_y(m))(u) &=& (f_y(m)\circ \mathscr{C}(v,x))(u)\\
&=&f_y(m)(uv)\\
&=&g(M(uv)(m))\\
&=& f_z(M(v)(m))(u).
\end{array}\vspace{-2pt}$$
This establishes our claim. Hence, the $R$-linear maps $f_y:M(y)\to I_x(y)$ with $y\in\mathscr{C}$ form a $\mathscr{C}$-linear morphism $f: M\to I_x$. Moreover, for any $m\in M(x)$, we have \vspace{-2pt} $$\Psi_{M, I_x}(f)(m)=\psi_x (f_x(m))=f_x(m)(1_x)=(g\circ M(1_x))(m)=g(m).\vspace{-1pt}$$
That is, $\Psi_{M, I_x}(f)=g$. So, $\Psi_{M, I_x}$ is surjective. Finally, a routine verification shows that 
$\Phi_{M, I_x}$ is natural in $M$. The proof of the proposition is completed.


\begin{Remark}\label{inj_Cmod}

Since the dual functor $D$ is exact, we deduce from Proposition \ref{Mor-mod-inj} that $I_x$ is an injective object in $\Mod \mathscr{C},$ for all $x\in \mathscr{C}$.
    
\end{Remark}


The following result was stated without any proof in \cite{Wat}. We include a short proof for the reader's convenience.

\begin{Prop}\label{Inj_C_mcat}

Let $\mathscr{C}$ be a small $R$-category. Then, $\Mod\mathscr{C}$ has enough injective objects.

\end{Prop}

\noindent{\it Poof.} Fix $M\in \Mod\mathscr{C}$. For $x\in \mathscr{C}$ and $\varphi\in D(M(x)),$ we write $x_\varphi=x$. Then, \vspace{.5pt} $I=\prod_{x\in \mathscr{C};\hspace{.5pt} \varphi\in D(M(x))}D\mathscr{C}(-,x_\varphi)$ is an injective object in $\Mod \mathscr{C};$ see (\ref{inj_Cmod}). For any $y\in \mathscr{C}$, we have an $R$-linear map \vspace{-2.5pt}
$$\psi_y\hspace{-1pt}: \hspace{-1pt} M(y)\to I(y)={\textstyle\prod}_{x\in \mathscr{C};\hspace{.5pt} \varphi\in D(M(x))} \, D\mathscr{C}(y,x_\varphi): m\mapsto (f_{x, \varphi,m})_{x\in\mathscr{C};\hspace{.5pt} \varphi\in D(M(x))}, \vspace{-2.5pt}
$$ where $f_{x,\varphi,m}\in D\mathscr{C}(y,x_\varphi)$ such that  $f_{x,\varphi,m}(u)=\varphi(M(u)(m))$, \vspace{.5pt} for $u\in \mathscr{C}(y,x_\varphi).$ Given a morphism $v:y\to z$ in $\mathscr{C}$, it is not hard to verify that the diagram \vspace{-4pt}
$$\hspace{30pt}\xymatrixcolsep{35pt}\xymatrixrowsep{22pt}\xymatrix{M(y)\ar[d]_{M(v)}\ar[r]^-{\psi_y} & \prod_{{x\in \mathscr{C};\hspace{.5pt} \varphi\in D(M(x))}}\mathscr{C}(y,x_\varphi) \ar@<-12ex>[d]^{\prod_{x\in \mathscr{C}; \hspace{.5pt} \varphi\in D(M(x))}D\mathscr{C}(x_\varphi, v)} \hspace{10pt} \\ M(z) \ar[r]^-{\psi_z} 
& \prod_{{x\in \mathscr{C}; \hspace{.5pt} \varphi\in D(M(x))}}D\mathscr{C}(z,x_\varphi)}\vspace{-4pt}$$ commutes. Assume that $\psi_y(m)=0$
for some $m\in M(y)$ with $y\in Q_0$. Then, for any $\varphi\in D(M(y))$, \vspace{1pt} we have $0=f_{y,\varphi,m}(1_y)=\varphi(M(1_y)(m))=\varphi(m)$. 
Since $I_R$ is a cogenerator for $\Mod R$, we have $m=0$. So, $\psi_y$ is a monomorphism. Therefore, the $R$-linear maps $\psi_y$ with $y\in \mathscr{C}$ form a $\mathscr{C}$-linear monomorphism $\psi: M\to I.$ The proof of the proposition is completed.

\subsection{\sc The Nakayama functor} Let $\mathscr{C}$ be a small $R$-category. We shall construct a Nakayama functor for $\Mod\mathscr{C}.$ Write ${\rm proj} \hspace{.5pt} \mathscr{C}$ and ${\rm inj}\mathscr{C}$ for the strictly full additive subcategories of $\Mod\mathscr{C}$ generated by the projective modules $P_x$ with $x\in \mathscr{C};$ see (\ref{proj_Cmod}), and by the injective modules $I_x$ with $x\in \mathscr C;$ see (\ref{inj_Cmod}), respectively.

\begin{Theo}\label{NK_functor}

Let $\mathscr{C}$ be a small $R$-category. Then, there exists a Nakayama functor
$\nu: {\rm proj}\mathscr{C}\to \Mod\mathscr{C},$ sending $P_x$ to $I_x$ for all objects $x\in \mathscr{C}.$
    
\end{Theo}

\noindent{\it Proof.} We begin with the full subcategory $\mathcal{P}$ of $\Mod \mathscr{C}$ generated by the modules $P_x$ with $x\in \mathscr{C}$. By the Yoneda Lemma, $\Hom_\mathscr{C}(P_x, P_y)=\{\mathscr{C}(u, -) \mid u\in \mathscr{C}(y,x)\}.$ Thus, setting $\nu P_x=I_x$  and $\nu \hspace{.5pt} \mathscr{C}(u, -)=D\mathscr{C}(-,u)$ for objects $x$ and morphisms $u$ in $\mathscr{C},$ we obtain an $R$-linear functor $\nu: \mathcal{P}\to \Mod\mathscr{C}.$ 

Fix $M\in \Mod\mathscr{C}$ and $x\in \mathscr{C}$. By Proposition \ref{Mor_proj_mod}, we have an $R$-linear isomorphism $\Phi_{P_x, M}: \Hom_{\mathscr{C}}(P_x, M) \to M(x),$ which is natural in $P_x$ and $M$. This induces an $R$-linear isomorphism $D\Phi_{P_x, M}: D(M(x)) \to D\Hom_{\mathscr{C}}(P_x, M),$ which is natural in $P_x$ and $M$. And by Proposition \ref{Mor-mod-inj}, we have an $R$-linear isomorphism 
\vspace{-2.5pt} $$\Psi_{M, \hspace{1pt} \nu P_x}: \Hom_{\mathscr{C}}(M, \nu P_x)\to D(M(x)): f\mapsto \theta_x\circ f_x,\vspace{-2.5pt}$$ which is natural in $M$. We claim that $\Psi_{M, \hspace{1pt} \nu P_x}$ is natural in $P_x$, or equivalently,
$\Psi_{M, \hspace{1pt} \nu P_y} \circ \Hom_{\hspace{.5pt}\mathscr{C}}(M, \hspace{1pt} \nu \hspace{.8pt}\mathscr{C}(u, -))= D(M(y)) \circ \Psi_{M, \hspace{1pt} \nu P_x},$ \vspace{1pt} for any $u: y\to x$ in $\mathscr{C}.$ 
Indeed, given $f\in \Hom_{\mathscr{C}}(M, \nu P_x)$ and $m\in M(y)$, we have \vspace{-2pt}
$$\hspace{-25pt} \begin{array}{rcl}
(\Psi_{M, \hspace{.5pt} \nu P_y}\circ \Hom_{\mathscr{C}}(M, \hspace{.5pt} \nu \hspace{1pt} \mathscr{C}(u, -)))(f)(m) \vspace{.5pt} 
&=& (\theta_y \circ (D\mathscr{C}(-,u)\circ f)_y)(m)\vspace{.5pt} \\
&=& \theta_y(D\mathscr{C}(y,u)(f_y(m)) \vspace{.5pt} \\
&=& (f_y(m)\circ \mathscr{C}(y,u))(1_y)\vspace{.5pt} \\
&=&f_y(m)(u).\end{array} \vspace{-3pt} $$
And since $f$ is $\mathscr{C}$-linear, we have $D\mathscr{C}(u, x) \circ f_y=f_x \circ M(u).$ This yields \vspace{-2pt}
$$\hspace{25pt} \begin{array}{rcl}
(D(M(u))\circ \Psi_{M, \hspace{.5pt} \nu P_x}(f))(m) 
&=&(\theta_x\circ D\mathscr{C}(u,x)\circ f_y)(m) \vspace{.5pt} \\
&=&\theta_x(D\mathscr{C}(u,x)(f_y(m))) \vspace{.5pt} \\
&=&(f_y(m) \circ \mathscr{C}(u,x))(1_x)\\
&=&f_y(m)(u).
\end{array} \hspace{10pt}\vspace{-3pt}$$
This establishes our claim. Composing $\Psi_{M, \hspace{.8pt} \nu P_x}$ and $D\Phi_{P_x, M}$ yields an $R$-linear isomorphism
$\Ta_{\hspace{-1pt} M, P_x}\hspace{-1.5pt}: \hspace{-1pt}\Hom_{\mathscr{C}}(M, \nu P_x) \hspace{-1pt}\to\hspace{-1pt} D\Hom_{\mathscr{C}}(P_x, M),$ which is natural in $P_x$ and $M$. Now, since every object in ${\rm proj}\mathscr{C}$ is isomorphic to a finite direct sum of the $P_x$ with $x\in \mathscr{C}$, we can extend, in a canonical way, the functor $\nu: \mathcal{P}\to \Mod\mathscr{C}$  to a functor $\nu: {\rm proj}\mathscr{C}\to \Mod \mathscr{C}.$ Moreover, the isomorphism 
$\Ta_{M, P_x}$ can also be extended to an $R$-linear isomorphism 
$\Ta_{M, P}: \Hom_{\mathscr{C}}(M, \nu P) \to D\Hom_{\mathscr{C}}(P, M),$ for 
 $P\in {\rm proj}\mathscr{C}$ and $M\in \Mod\mathscr{C}$, and it is natural in $P$ and $M$. This shows that $\nu:  {\rm proj}\mathscr{C}\to \Mod \mathscr{C}$ is a Nakayama functor. The proof of the theorem is completed.

\vspace{-2pt}

\subsection{\sc Finitely presented and finitely copresented modules} 
Let $\mathscr{C}$ be a small $R$-category. Fix a module $M$ in $\Mod\mathscr{C}$. We shall say that $M$ is {\it finitely generated} if there exists an epimorphism $f: P\to M$ with $P\in {\rm proj}\hspace{.5pt}\mathscr{C};$ and {\it finitely presented} if it admits a projective presentation over ${\rm proj} \hspace{.5pt} \mathscr{C},$ or equivalently, there exists an epimorphism $f: \hspace{-1.5pt} P \hspace{-1.5pt} \to \hspace{-1.5pt}  M$ with  $P$ in ${\rm proj}\hspace{.5pt}\mathscr{C}$ and ${\rm Ker}(f)$ finitely generated. Dually, $M$ is called {\it finitely cogenerated} if there exists a monomorphism $g: M\to I$ with $I\in {\rm inj} \hspace{.5pt}\mathscr{C};$ and {\it finitely copresented} if it admits an injective copresentation over ${\rm inj} \hspace{.5pt} \mathscr{C},$
or equivalently, there exists a monomorphism $g:  M\to I$ with $I$ in ${\rm inj}\hspace{.5pt}\mathscr{C}$ and ${\rm Coker}(g)$ finitely cogenerated. We denote by $\mmod^+\mathscr{C}$ and $\mmod^-\mathscr{C}$ the full subcategories of finitely presented modules and finitely copresented modules in $\Mod \mathscr{C},$ respectively. By Propositions \ref{subcat_pres} and \ref{subcat_copres}, $\mmod^+\mathscr{C}$ and $\mmod^-\mathscr{C}$ are exact $R$-categories. It is interesting to know when they are abelian categories. 

\vspace{-0pt}

\begin{Prop}\label{fpres_AB}

Let $\mathscr{C}$ be a small $R$-category. Then 

\begin{enumerate}[$(1)$]

\vspace{-2pt}

\item ${\rm mod}^{+}\mathscr{C}$ is abelian if and only if the finitely generated submodules of the modules in ${\rm proj}\hspace{.8pt}\mathscr{C}$ are finitely presented$\hspace{.8pt};$

\vspace{.5pt}

\item ${\rm mod}^{-}\mathscr{C}$ is abelian if and only if the finitely cogenerated quotient-modules of the modules in ${\rm inj}\hspace{.8pt}\mathscr{C}$ are finitely copresented.
    
\end{enumerate}\end{Prop}

\vspace{-1.5pt}

\noindent{\it Proof.} We shall only prove Statement (1). Assume that ${\rm mod}^{+}\mathscr{C}$ is abelian. Let $M$ be a finitely generated submodule of a module $P_0$ in ${\rm proj}\hspace{.8pt} \mathscr{C}.$ Then, we have an epimorphism $f: P_1\to M$ with $P_1\in {\rm proj}\hspace{.8pt} \mathscr{C}$, and hence, a morphism $jf: P_1\to P_0$ in ${\rm mod}^+ \mathscr{C}$, where $j: M\to P_0$ is the inclusion map. Note that $f$ and $jf$ have the same kernel $K$, which lies in ${\rm mod}^{+} \mathscr{C}.$ Thus, $M$ is finitely presented. 

Assume that the finitely generated submodules of modules in ${\rm proj}\hspace{.8pt}\mathscr{C}$ are finitely presented. Consider a morphism $f: P_1\to P_0$ in ${\rm proj}\hspace{.5pt} \mathscr{C}$ and its kernel $j: K\to P_1$ in $\Mod \mathscr{C}.$ By the assumption, there exists an epimorphism $f': P_1'\to {\rm Im}(f)$ having a kernel $q: K'\to P_1'$, where $P_1'\in {\rm proj}\hspace{.8pt} \mathscr{C}$ and $K'$ is finitely generated. By Schanule's Lemma, $K'\oplus P_1\cong K\oplus P_1'.$ \vspace{-3pt} Thus, $K$ is finitely generated with an epimorphism $g: P_2\to K,$ where $P_2\in {\rm proj}\hspace{.8pt}\mathscr{C}$. So, $\xymatrixcolsep{20pt}\xymatrix{\hspace{-4pt} P_2\ar[r]^{jg} & P_1\ar[r]^f & P_0}\hspace{-1pt}\vspace{-1pt}$ is an exact sequence in ${\rm proj}\hspace{.8pt} \mathscr{C}$. Thus, ${\rm mod}^{+}\mathscr{C}$ is closed under kernels; see \cite[(2.1)]{Aus1}. \vspace{.5pt} Being also closed under cokernels; see (\ref{subcat_pres}),  ${\rm mod}^{+}\mathscr{C}$ is abelian. 
%
The proof of the proposition is completed.



\vspace{-3pt}

\subsection{\sc The existence thereom} In this subsection, we need to assume that $\mathscr{C}$ is endo-local and Hom-reflexive. In this setting, we have the following statement. 

\begin{Prop}\label{C-proj}

Let $\mathscr{C}$ be a small endo-local Hom-reflexive $R$-category. Then

\begin{enumerate}[$(1)$]

\vspace{-2.5pt} 

\item ${\rm proj}\mathscr{C}\hspace{-.5pt}$ and ${\rm inj}\mathscr{C}\hspace{-.5pt}$ are Hom-reflexive Krull-Schmidt subcategories of $\mmod \mathscr{C};$ 



\vspace{.5pt} 

\item the functor $\mk D: \Mod \mathscr{C}^\circ \hspace{-1.5pt}  \to \hspace{-1.5pt} \Mod\mathscr{C}$ restricts to a duality $\mf{D}: {\rm mod}^{+}\mathscr{C}^\circ \hspace{-1pt} \to \hspace{-1pt} {\rm mod}^{-}\mathscr{C}\hspace{-1.5pt}.$ 
    
\end{enumerate}

\end{Prop} 

\noindent{\it Proof.} Let $x, y\in \mathscr{C}$. Since $P_x(y)=\mathscr{C}(x, y),$ 
${\rm proj}\mathscr{C}\hspace{-.5pt}\subseteq \mmod \mathscr{C}\hspace{-2pt}.$ By Proposition \ref{Mor_proj_mod}, \vspace{.5pt} there exists an $R$-linear isomorphism $\varphi_{x,y}: \Hom_{\mathscr{C}}(P_x, P_y)\to \mathscr{C}(y,x)$, natural in $x$ and $y$. So, 
${\rm proj}\mathscr{C}$ is Hom-reflexive. Since $\varphi_{x,x}: \End_{\mathscr{C}}(P_x)\to \mathscr{C}(x,x)$ is an algebra anti-isomorphism, ${\rm proj}\mathscr{C}$ is Krull-Schmidt. 
Since $\mathscr{C}^\circ$ is Hom-reflexive and Krull-Schmidt, ${\rm proj}\mathscr{C}^\circ$ is Hom-reflexive, Krull-Schmidt and contained in $\mmod \mathscr{C}^\circ$. Then, the duality $\mk D: \mmod \mathscr{C}^\circ \to \mmod \mathscr{C}$ restricts a duality $\mk D: {\rm proj}\mathscr{C}^\circ \to {\rm inj}\mathscr{C};$ see (\ref{LRMod}) Therefore, ${\rm inj}\mathscr{C}\hspace{-.5pt}$ is Hom-reflexive, Krull-Schmidt and contained in $\mmod \mathscr{C}.$ 

Since $\mmod \mathscr{C}^\circ$ and $\mmod \mathscr{C}$ are abelian, ${\rm mod}^{+}\mathscr{C}^\circ \subseteq \mmod \mathscr{C}^\circ$ and ${\rm mod}^{-}\mathscr{C}\subseteq \mmod \mathscr{C}.$ For any $M\in {\rm mod}^{+}\mathscr{C}^\circ$ with a projective presentation 
$\xymatrixcolsep{18pt}\xymatrix{\hspace{-2pt}P_1 \ar[r]& P_0 \ar[r] &M \ar[r] &0}\vspace{-2.5pt}\hspace{-2pt}$ over ${\rm proj}\hspace{.5pt}\mathscr{C}^\circ\hspace{-2pt},$ we have an injective copresentation $\xymatrixcolsep{18pt}\xymatrix{\hspace{-2pt} 0\ar[r] &\mk DM\ar[r]&\mk DP_0\ar[r]&\mk DP_1}\vspace{-1.5pt}\hspace{-2pt}$ over ${\rm inj}\hspace{.5pt}\mathscr{C}\hspace{-1.5pt}.$ Hence, $\mk DM\in {\rm mod}^{-}\mathscr{C}$. So, the duality $\mk D: \mmod \mathscr{C}^\circ \to \mmod \mathscr{C}$ restricts to a duality $\mf{D}: {\rm mod}^{+}\mathscr{C}^\circ \hspace{-1pt} \to \hspace{-1pt} {\rm mod}^{-}\mathscr{C}\hspace{-1.5pt}.$  The proof of the proposition is completed.

\begin{Remark} Let $\mathscr{C}$ be a small $R$-category, where $R$ is a commutative artinian ring. By Proposition \ref{C-proj}(2), we see that $\mathcal{C}$ is a dualizing $R$-variety as defined in \cite[Section 2]{AIR} if and only if $\mathscr{C}$ is endo-local and Hom-finite with $\mmod^+\mathscr{C}=\mmod^-\mathscr{C}.$
    
\end{Remark}

\vspace{3pt}

We are ready to obtain the main result of this section.

\begin{Theo}\label{C_mod_ass}

Let $\mathscr{C}$ be a small endo-local Hom-reflexive $R$-category. Consider a strongly indecomposable module $M$ in $\Mod\mathscr{C}.$

\begin{enumerate}[$(1)$]

\vspace{-1pt}

\item If $M\in \mmod^+\mathscr{C}$ is not projective, \vspace{-1.5pt} then there exists an almost split sequence  $\xymatrixcolsep{20pt}\xymatrix{\hspace{-1pt} 0\ar[r] &\tau M\ar[r] & E \ar[r] &M \ar[r] &0 \hspace{-2pt}}$ in $\Mod\mathscr{C}$, where $\tau M\in \mmod^-\mathscr{C}.$

\item If $M\in \mmod^-\mathscr{C}$ is not injective, \vspace{-1.5pt} then there exists an almost split sequence $\xymatrixcolsep{20pt}\xymatrix{\hspace{-1pt} 0\ar[r] & M\ar[r] & E \ar[r] &\tau^-\!M \ar[r] &0 \hspace{-2pt}}$  in $\Mod\mathscr{C}$ where $\tau^- \! M\in \mmod^+\mathscr{C}.$

\end{enumerate} \end{Theo}

\vspace{-1pt}

\noindent{\it Proof.} By Theorem \ref{NK_functor}, we have a Nakayama functor $\nu: {\rm proj}\mathscr{C}\to \Mod\mathscr{C},$ where ${\rm proj}\mathscr{C}$ is Hom-reflexive and Krull-Schmidt; see (\ref{C-proj}). Since $\Mod\mathscr{C}$ has enough projective and injective objects; see (\ref{Proj_C_mcat}) and (\ref{Inj_C_mcat}), the result follows immediately from Theorem \ref{AR-seq}. The proof of the theorem is completed.

\begin{Remark}\label{ass_cln}

Let $\mathscr{C}$ be a small endo-local Hom-noetherian $R$-category  
with $R$ being complete noetherian local. 

\vspace{1pt}

\noindent (1) In Theorem \ref{C_mod_ass}, it suffices to assume that $M$ is indecomposable; see (\ref{fpres_cat_KS}).

\vspace{1pt}

\noindent (2) Theorem \ref{C_mod_ass}(1) is stated 
as Theorem 6 in \cite{Aus2}. However, we cannot find any proof in the existing literature, and we do not see how to establish an Auslander-Reiten duality for $\Mod\mathscr{C}$ by the classical approach as used in \cite[(7.4)]{AIR} or \cite[(I.3.4)]{A2}.

\end{Remark}

\vspace{1pt}

\section{Application to representations of quivers with relations}

The objective of this section is mainly to apply our previous results to study the existence of almost split sequences over algebras defined by quivers with relations. In the local semiperfect case; see (\ref{lspa}), a general existence theorem in the category of all modules is a special case for small endo-local Hom-reflexive categories; see (\ref{ass_lfd}). So, we shall focus on the subcategories of finitely presented modules, finitely copresented modules and finite dimensional modules. For locally semiprimary algebras given by locally finite quivers; see (\ref{lspri_def}), we show that the almost split sequences in the first two subcategories are almost split in the category of all modules; see (\ref{ass_fpres}) and (\ref{ass_fcp}). More importantly, we shall study when these subcategories have almost split sequences; see (\ref{ass_fp_cat}), (\ref{ass_fcp_cat}) and (\ref{ass_fdm_cat}). 

\vspace{1pt}

Throughout this section, $k$ denotes a field. Note that the dual functor for $\Mod k$ is $D=\Hom_k(-, k),$ and a module in $\Mod k$ is reflexive if and only if it is finite dimensional. In particular, a $k$-category is Hom-reflexive if and only if it is Hom-finite. In this section, a $k$-algebra has enough idempotents but not necessarily an identity, in which an ideal is always two-sided. 

\vspace{-2pt}

\subsection{\sc Quivers} Let $Q=(Q_0, Q_1)$ be a quiver, where $Q_0$ is the set of vertices and $Q_1$ is the set of arrows. Given $\alpha: x\to y$ in $Q_1$, we call $x$ the {\it starting point} and $y$ the {\it ending point}, and write $s(\alpha)=x$ and $e(\alpha)=y$. For each $x\in Q_0$, one associates a {\it trivial path} $\varepsilon_x$ with $s(\varepsilon_x)=e(\varepsilon_x)=x$. A {\it path} of length $n \hspace{.5pt} (>0)$ is a sequence $\rho=\alpha_n \cdots \alpha_1$, where $\alpha_i\in Q_1$ such that $s(\alpha_{i+1})=e(\alpha_i)$, for $1\le i< n$. Two paths are called {\it comparable} if one is the subpath of the other, and a set of paths is called {\it comparable} if every two paths in the set are comparable. For $x\in Q_0$ and $n\ge 0$, we denote by $Q_n(x,-)$ and $Q_n(-,x)$ the sets of paths of length $n$ starting with $x$, and those ending with $x$, respectively. We say that $Q$ is {\it left locally finite} if $Q_1(x, -)$ is finite for all $x\in Q_0;$ {\it right locally finite} if $Q_1(-,x)$ is finite for all $x\in Q_0,$ and {\it locally finite} if $Q$ is left and right locally finite. Further, $Q$ is called {\it interval-finite} if the set of paths from $x$ to $y$ is finite, for all $x, y\in Q_0$; and {\it strongly locally finite} if $Q$ is locally finite and interval-finite.

The {\it opposite} quiver $Q^{\rm o}$ of $Q$ is defined in such a way that $(Q^{\rm o})_0=Q_0$ and $(Q^{\rm o})_1=\{\alpha^{\rm o}: y\to x \mid \alpha: x\to y \in Q_1\}$. A non-trivial path $\rho=\alpha_n\cdots \alpha_1$ in $Q$, where $\alpha_i\in Q_1$, induces a non-trivial path $\rho^{\rm o}=\alpha_1^{\rm o} \cdots \alpha_n^{\rm o}$ in $Q^{\rm o}$. However, we shall identify the trivial path in $Q^{\rm o}$ at a vertex $x$ with that in $Q$.

\vspace{-1pt}

\subsection{\sc Algebras given by quivers with relations} Let $Q$ be a quiver. Write $kQ$ for the path algebra of $Q$ over $k$. An ideal $\mk I$ in $kQ$ is a {\it relation ideal} if $\mk I\subseteq (kQ^+)^2$, where $kQ^+$ is the ideal in $kQ$ generated by the arrows in $Q$. In this case, $(Q, \mk{I})$ is called a {\it bound quiver}, and a path $\rho$ in $Q$ is called {\it nonzero in $\La$}, or simply {\it nonzero}, if $\rho\not\in \mk I.$ Given $\rho=\sum \lambda_i \rho_i\in kQ$, where $\lambda_i\in k$ and $\rho_i$ is a path in $Q$, we write $\rho^\circ=\sum \lambda_i \rho_i^\circ\in kQ^\circ$. In this way, an ideal $\mk I$ in $kQ$ is a relation ideal if and only if $\mk I^\circ=\{\rho^\circ \mid \rho \in \mk I\}$ is a relation ideal in $kQ^\circ.$

\vspace{0pt}

Consider now a $k$-algebra $\La=kQ/\mk{I},$ where $Q$ is a quiver and $\mk I$ is a relation ideal in $kQ$. 
Write $\bar{\hspace{-1pt}\gamma}=\gamma+\mk I\in \La$ for $\gamma\in kQ$. Then $\La$ has a complete orthogonal set of idempotents $\{e_x \mid x\in Q_0\},$ where $e_x= \bar\varepsilon_x$. Set $J_\mathit{\hspace{-1pt}\Lambda}=\{\bar \rho \mid \rho \in kQ^+\}:=J.$ Note that $kQ^\circ/ \mk I^\circ$ is the opposite algebra of $\La$. Write $\bar \rho^\circ=\rho^\circ + \mk{I}^\circ$ for $\rho\in kQ$ but $e_x=\varepsilon_x+\mk I^\circ$ for $x\in Q_0$. In particular, $J_{\hspace{-1pt}\mathit\Lambda^\circ}=\{\bar \rho^\circ \mid \rho \in kQ^+\}=J^\circ$. 

\vspace{1pt}

A left $\La$-module $M$ is called {\it unitary} if $M=\sum_{x\in Q_0} e_xM$. We denote by $\Mod\La$ the category of all unitary left $\La$-modules. Observe that the $k$-algebra $\La$ can also be regarded as a small $k$-category, whose objects are the vertices of $Q$ and whose morphism spaces are $\La(x, y)=e_y\La e_x$ with $x, y\in Q_0$. In this viewpoint, a $k$-linear representation of $(Q, \mk{I})$ is a covariant functor from $\La$ to $\Mod k.$ It is well known that $\Mod \La$ is equivalent to the category ${\rm Rep}(Q, \mk I)$ of all $k$-linear representations of $(Q, \mk{I}),$ in such a way that a module $M$ corresponds to a representation $\tilde{M}$ defined by
$\tilde{M}(x)=e_xM$ and $\tilde{M}(u): \tilde{M}(x)\to \tilde{M}(y)$ \vspace{-1pt} being the left multiplication by $u$, where $x, y\in Q_0$ and $u \in \La(x, y)$.
Identifying $M$ with $\tilde{M}$, we shall freely use the terminology, the notation and the results stated in Section 4 for $\La$-modules. In particular, for each  $x\in Q_0$, we have a projective module $P_x=\La e_x=\La(x,-);$ see (\ref{proj_Cmod}) and a simple module $S_x=\La e_x/J e_x$ in $\ModLa.$ And we have a projective module $P_x^\circ=\La^\circ e_x=\La(-, x)$ and a simple module $S_x^\circ=\La^\circ e_x/ J^\circ e_x$ in $\ModLa^\circ.$ 


For a module $M$ in $\ModLa^\circ$, we have a module
$\mf{D}M = \oplus_{x\in Q_0}\, D(e_xM)$ in $\ModLa$ such, for all $\varphi\in D(e_xM)$ and $u\in e_y\La e_x$, that $u \hspace{.5pt} \varphi\in D(e_yM)$, which is defined by $(u \varphi)(m)=\varphi(u^{\rm o} m),$ for $m\in e_yM.$ For a morphism $f: M\to N$ in $\ModLa^\circ,$ we have a morphism $\mf{D}f: \mf{D}N\to \mf{D}M$ in $\ModLa$ such, for $\varphi\in D(e_xN)$ with $x\in Q_0$, that $(\mk D f)(\varphi)= \varphi \circ f_x$, where $f_x: e_xM\to e_xN$ is the restriction of $f$. This yields a contravariant functor $\mf{D} : \ModLa^{\rm o} \to \ModLa.$ For each $x\in Q_0,$ we see that $\m{D}(S_x^\circ)\cong S_x$ and $I_x=\mk{D}(P_x^\circ)=\mk D(-,x)$, an injective module in $\ModLa;$ see (\ref{inj_Cmod}).

\vspace{1pt}

A module $M$ in $\ModLa$ is \vspace{-1pt} called {\it locally finite dimensional} if ${\rm dim}_k(e_xM)<\infty$ for all $x\in Q_0.$ We denote by ${\rm mod}\La$ and ${\rm mod}^{\hspace{.4pt}b\hspace{-2.5pt}}\La$ the full subcategories of locally finite dimensional modules and finite dimensional modules in $\Mod \La$, respectively. As a special case of Proposition \ref{LRMod}, we obtain the following statement. 


\begin{Prop}\label{duality_lfdm}

Let $\La=kQ/\mk I$, where $Q$ is a quiver and $\mk I$ is an ideal in $kQ$. Then the contravariant functor $\mk D: \ModLa^\circ \to \ModLa$ restricts to dualities $\mk D: \mmod \La^\circ \to \mmod \La$ and $\mk D: {\rm mod}^{\hspace{.4pt}b\hspace{-2.5pt}}\La^\circ \to {\rm mod}^{\hspace{.4pt}b\hspace{-2.5pt}}\La$.

\end{Prop}

\vspace{-2pt}

\subsection{\sc Locally noetherian algebras} Let $\La=kQ/\mk I$, where $Q$ is a quiver and $\mk I$ is an ideal in $kQ$. 
Recall that $\La$ is {\it locally left noetherian} if $\La e_x$ is noetherian for all $x\in Q_0$, {\it locally right noetherian} if $e_x\La$ is noetherian for all $x\in Q_0$, and {\it locally noetherian} if $\La$ is locally left and right noetherian. Moreover, $\La$ is {\it locally left bounded} if $\La e_x$ is finite dimensional for all $x\in Q_0,$ {\it locally right bounded} if $e_x\La$ is finite dimensional for all $x\in Q_0,$ and {\it locally bounded} if $\La$ is locally left and right bounded. Locally (left, right) bounded algebras are clearly locally (left, right) noetherian, and so are 
the special multiserial algebras; see \cite[(1.3)]{LL25}. 
To provide more classes of locally noetherian algebras, we introduce the following notion.

\vspace{-1pt}

\begin{Defn}\label{mserial}

Let $\La=kQ/\mk{I},$ where $Q$ is a quiver and $\mk I$ is a relation ideal in $kQ$. We shall say that $\La$ is

\begin{enumerate}[$(1)$]

\vspace{-1.5pt}
    
\item {\it locally left eventually multiserial} if $Q$ is left locally finite, and for each $x\in Q_0$, there exists an integer $n_x (\ge 0)$ such, for any $\rho \in Q_{n_x}(x, -)$, that the set of nonzero paths in $Q$ starting with $\rho$ is comparable$\,;$

\vspace{.5pt}

\item {\it locally right eventually  multiserial} if $Q$ is right locally finite, and for any $x\in Q_0$, there exists an integer $n_x (\ge 0)$ such, for any $\rho\in Q_{n_x}(-, x)$, that the set of nonzero paths in $Q$ ending with $\rho$ is comparable$\,;$

\item {\it locally eventually  multiserial} if $\La$ is locally left and right eventually multiserial.

\end{enumerate} \end{Defn}

\vspace{.5pt}

Note that special multiserial algebras are locally left and right eventually multiserial, while locally bounded algebras are not necessarily so. 

\begin{Prop}\label{mul-noe}

Let $\La=kQ/\mk I,$ where $Q$ is a quiver and $\mk I$ is a relation ideal in $kQ$. If $\La$ is locally left or right eventually multiserial, then it is locally left or right noetherian, respectively.

\end{Prop}


\noindent{\it Proof.} We shall consider only the case where $\La$ is locally left eventually multiserial. Let $x$ be an vertex in $Q$ with an integer $n_x$ as stated in Definition \ref{mserial}(1). Fix $\rho\in Q_{n_x}(x, -)$. Then \vspace{.5pt} $\La \bar \rho=\sum_{\hspace{.5pt}\eta\in \mathit\Omega(\rho)} k \bar\eta,$ where $\Oa(\rho)$ denotes the set of nonzero paths in $Q$ starting with $\rho$. If $\Oa(\rho)$ is finite, then $\La\bar\rho$ is noetherian. Otherwise, since $\Oa(\rho)$ is comparable, there exist arrows $\alpha_i: x_i\to x_{i+1}$ in $Q$ with $i\ge 1$ and $x_1=x,$ such that 
$\Oa(\rho)=\{\rho, \alpha_1\rho, \ldots, \alpha_i \cdots \alpha_1\rho, \ldots \}.$ Set 
$u_0=\bar \rho$ and $u_i=  \bar\alpha_i \cdots \bar \alpha_1 \bar \rho$ for $i\ge 1$. Given $0\ne u\in \La\bar\rho$, there exists a minimal $d(u)\ge 0$ such that $u=\sum_{i=0}^{d(u)} \lambda_i u_i$, where $\lambda_i\in k$ with $\lambda_{d(u)}\ne 0$. If $u, v\in \La\bar\rho$ are nonzero, then $v=qu + w$, where $q\in \La$ and $w\in \La \bar \rho$ with $w=0$ or $d(w)<d(u).$ Thus, every nonzero submodule of $\La \bar \rho$ is generated by some $u\in \La\bar\rho$. Thus, $\La \bar\rho$ is noetherian. Since $Q$ is left locally finite, $J^{n_x}e_x=\sum_{\rho\in Q_{n_x}(x, -)} \hspace{-1pt} \La \bar \rho$ \vspace{1pt} with $Q_{n_x}(x, -)$ finite is noetherian, and $\La e_x/J^{n_x}e_x$ is finite dimensional. So, $\La e_x$ is noetherian. The proof of the proposition is completed.

\begin{Exam}\label{trun_mser_ex} 

Consider $\La=kQ/\mk I$, where $Q$ is the quiver \vspace{-2.5pt} 
$$\xymatrixrowsep{20pt}\xymatrix{
\ar@{.}[r] & a_3\ar[r]^{\alpha_3} & a_2\ar[r]^{\alpha_2} & a_1 \ar[r]^{\alpha_{1}}  & a_0\ar@<.5ex>[r] \ar@<-.8ex>[r] \ar[d] & d_1 \ar[r] \ar[d]  & d_2 \ar[d] \ar[r] & d_3 \ar[d] \ar@{.}[r] & \\
& b_3\ar[u]_{\beta_3} & b_2\ar[u]_{\beta_2}& b_1\ar[u]_{\beta_1}&c_0 &c_1&c_2& c_3&} \vspace{-2pt}$$ and $\mk{I}=\langle \alpha_i\beta_i \mid i\ge 1\rangle$. Then $\La$ is locally right eventually multiserial. By Proposition \ref{mul-noe}, it is locally right noetherian. On the other hand, $\La$ is neither locally left eventually multiserial nor locally left noetherian.
\end{Exam}

\vspace{-5pt}

\subsection{\sc Locally semiperfect algebras} A ring $\Sa$ with an identity is semiperfect if and only if it has a complete orthogonal set of idempotents $\{e_1, \ldots, e_n\}$ such that $e_i \Sa e_i$ is local; see \cite[(27.6)]{AnF}. This inspires the following definition; see \cite[(1.4.3)]{Gup}.

\begin{Defn}\label{lspa}

Let $\La=kQ/\mk{I}$, where $Q$ is a quiver and $\mk{I}$ is a relation ideal in $kQ$. One says that $\La$ is {\it locally semiperfect} provided, for all $x, y\in Q_0$, that 

\begin{enumerate}[$(1)$]

\vspace{-3.5pt}

\item $e_x\La e_x$ is a local $k$-algebra;

\item $e_y \La e_x$ is a finite dimensional $k$-vector space. 

\end{enumerate}\end{Defn}

\begin{Remark} Let $\La=kQ/\mk{I}$ be a locally semiperfect algebra. Regarded as a $k$-category, $\La$ is small endo-local and Hom-finite.

\end{Remark}

\begin{Exam}\label{spa_ex} 

(1) Let $\La$ be the $k$-algebra defined by a single loop $\alpha$ at a vertex with relation
$\alpha^2=\alpha^3$. Then $\La$ is locally bounded but not locally semiperfect. 

\vspace{1pt}

\noindent (2) Consier the $k$-algebra $\La=kQ/\mk I$, where $Q$ is the quiver \vspace{-5pt}
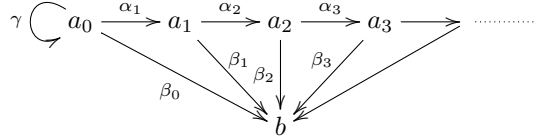
$$\hspace{70pt}\xymatrixcolsep{22pt}\xymatrixrowsep{28pt}\xymatrix{a_0\ar[r]^{\alpha_1} \ar[drr]_{\beta_0} \ar@(ul,dl)_\gamma & a_1 \ar[r]^{\alpha_2} \ar[dr]^{\hspace{-4pt}\beta_1} & a_2\ar[r]^{\alpha_3}\ar[d]_{\beta_2\hspace{-1pt}} &a_3\ar[dl]_{\beta_3\hspace{-4pt}} \ar[r]  &  \ar[dll] \ar@{.}[r] &\\
&& b &&&&&&&}\vspace{-7pt}$$
and $\mk{I}=\langle \gamma^2, \beta_1\alpha_1-\beta_i\alpha_i\cdots \alpha_2\alpha_1 \mid i\ge 2 \hspace{.5pt} \rangle.$ Then $\La$ is locally semiperfect with $Q$ being left locally finite but not right locally finite. 

\end{Exam}

\vspace{2pt}

The following result collects some basic properties of locally semiperfect algebras.

\begin{Prop}\label{lfd_rad}

Let $\La=kQ/\mk{I}$ be locally semiperfect, where $Q$ is a quiver and $\mk I$ is a relation ideal in $kQ$. 

\begin{enumerate}[$(1)$]

\vspace{-3pt}

\item If $x\in Q_0$, then $Je_x$ is the unique maximal submodule of $P_x$.

\item If $x\in Q_0$, then $\soc \hspace{.5pt} I_x$ is isomorphic to $S_x$ and essential in $I_x.$

\item The non-isomorphic simple modules in $\Mod \La$ are the
$S_x$ with $x\in Q_0$.

\item If $M\in \Mod\La,$ then 
$\soc M=\{ m\in M \mid J m=0\}.$

\end{enumerate} \end{Prop}

\vspace{-1pt}

\noindent{\it Proof.} (1) Let $x\in Q_0$. Clearly, $Je_x$ is a maximal submodule of $P_x$. For the uniqueness, it suffices to show that $e_x-\bar \rho$ is invertible for $\rho\in \varepsilon_x(kQ^+) \varepsilon_x$. Assume that this is not the case. Since $e_x\La e_x$ is local, $\bar \rho$ is invertible for some $\rho\in \varepsilon_x(kQ^+) \varepsilon_x$. Then, $e_x=\bar \rho \bar \eta$ for some $\eta \in \varepsilon_x(kQ) \varepsilon_x$. So $e_x-\rho\eta\in \mk I,$ contrary to $\mk I \subseteq (kQ^+)^2.$ 

(2) Let $x\in Q_0$. Since $\La^\circ$ is locally semiprimary, by Statement (1), the canonical projection $p: \La^\circ e_x\to S_x^\circ$ is a projective cover. So, $p$ is superfluous; see \cite[(3.4)]{Krau}. Since $\La^\circ e_x \in\mmod \La^\circ$, by Proposition \ref{duality_lfdm}, $\mk Dp: S_x\to I_x$ is an essential monomorphism. Thus, ${\rm Im}(\mk Dp)$ is an essential submodule of $I_x$. Since ${\rm Im}(\mk Dp)\cong S_x$, we see that $\soc I_x={\rm Im}(\mk Dp)$.

(3) Let $S$ be a simple module in $\ModLa.$ Then, there exists $0\ne m\in e_xS$ for some $x\in Q_0$. If $(Je_x)m=S$, then $m=\bar \rho m$ for some $\bar\rho \in e_x J e_x.$ Since $e_x-\bar\rho$ is invertibe, we have $m=0$, absurd. Thus, $(Je_x) m=0,$ and hence, $S\cong \La e_x/ Je_x.$

(4) Given $m\in \soc M$, by Statement (3), $Jm=0$. Suppose that $Jm=0$ for some $0\ne m\in M.$ Then, $m=e_{x_1}m+\cdots +e_{x_r}m,$ where $e_{x_i}m\ne 0$ and $x_i\in Q$ are pairwise distinct. Since $Je_{x_i}(e_{x_i}m) = (Je_{x_i})m \subseteq Jm=0$, we deduce that $\La(e_{x_i}m_i)\cong S_{x_i}$, for $i=1, \ldots, r$. Thus, $m\in \soc M$. The proof of the proposition is completed.

\begin{Remark} By Proposition \ref{lfd_rad}(1), the condition (1) in Definition \ref{lspa} can be replaced by the condition that every oriented cycle in $Q$ is nilpotent in $\La.$

\end{Remark}


\begin{Lemma}\label{fd_soc}

Let $\La=kQ/\mk I$ be locally semiperfect, where $Q$ is a quiver and $\mk I$ is a relation ideal in $kQ$. Then a module $M$ in $\Mod \La$ is finitely cogenerated if and only if $\soc M$ is finite dimensional and essential in $M\hspace{.4pt}.$

\end{Lemma}

\noindent{\it Proof.} Assume that $M$ is finitely cogenerated with an injective envelope $g: M\to I,$ where $I\in {\rm inj}\La$. By Proposition \ref{lfd_rad}(2), $\soc I\in {\rm mod}^{\hspace{.4pt} b\hspace{-2.5pt}}\La$ is essential in $I$. Since $g(\soc M)\subseteq \soc I$, we have $\soc M\in {\rm mod}^{\hspace{.4pt} b\hspace{-2.5pt}}\La$ and $g(\soc M)=\soc I.$ Thus, $\soc M$ is essential in $M.$ Conversely, assume that $\soc M\in {\rm mod}^{\hspace{.4pt} b\hspace{-2.5pt}}\La$ is essential in $M.$ Then, $\soc M=km_1\oplus \cdots \oplus km_s$, where $0\ne m_i\in e_{x_i}(\soc M)$ with $x_i\in Q_0$. By Proposition \ref{lfd_rad}(4), $km_i\cong S_{x_i}$, and in view of Proposition \ref{lfd_rad}(2), we have a monomorphism $h_i: km_i\to I_{x_i}$ with ${\rm Im}(h_i)=\soc I_{x_i}$, for $i=1, \ldots, s.$ This yields a monomorphism $h: \soc M\to I=\oplus_{i=1}^s I_{x_i}$ with $h(\soc M)=\soc I.$ Thus, $h$ extends to a monomorphism $g: M\to I$. The proof of the lemma is completed.

\vspace{-2pt}

\subsection{\sc Locally semiprimary algebras} For our later study on almost split sequences in ${\rm mod}^{+\hspace{-3pt}}\La$ and  ${\rm mod}^{-\hspace{-3pt}}\La$, we need to impose some stronger conditions on $\La$. Recall that a ring $\Sa$ with an identity is called {\it semiprimary} if $\rad \Sa$ is nilpotent and $\Sa/\rad \Sa$ is semisimple; see \cite[(4.15)]{Lam}. This motivates the following definition.

\begin{Defn}\label{lspri_def}

Let $\La=kQ/\mk{I}$, where $Q$ is a quiver and $\mk{I}$ is a relation ideal in $kQ$. We say that $\La$ is {\it locally semiprimary} provided, for any $x, y\in Q_0$, that $Q$ contains at most finitely many nonzero paths from $x$ to $y$.

\end{Defn}

\begin{Remark}\label{lsprm_lsp} Clearly, locally semiprimary algebras are locally semiperfect.

\end{Remark}

\begin{Exam} \label{lspri_ex} (1) If $Q$ is an interval-finite quiver, then $kQ$ is locally semiprimary.

\vspace{2pt}

\noindent (2) Consider $\La=kQ/\mk I$, where $Q$ is the quiver \vspace{-2pt}
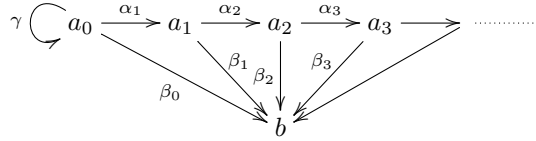
$$\hspace{70pt}\xymatrixcolsep{22pt}\xymatrixrowsep{28pt}\xymatrix{a_0\ar[r]^{\alpha_1} \ar[drr]_{\beta_0} \ar@(ul,dl)_\gamma & a_1 \ar[r]^{\alpha_2} \ar[dr]^{\hspace{-4pt}\beta_1} & a_2\ar[r]^{\alpha_3}\ar[d]_{\beta_2\hspace{-1pt}} &a_3\ar[dl]_{\beta_3\hspace{-4pt}} \ar[r]  &  \ar[dll] \ar@{.}[r] & \\
&& b &&&&&&&}\vspace{-4pt}$$
and $\mk{I}=\langle \gamma^2,  \beta_n\alpha_n\cdots \alpha_i\mid i\ge 1; n\ge 2i\hspace{.5pt}\rangle.$ Then $\La$ is locally semiprimary with $Q$ not locally finite. In contrast, the algebra in Example \ref{spa_ex}(2) defined by the same quiver is not locally semiprimary.
    
\end{Exam}
 
\vspace{1pt}

We shall need the following statement for describing almost split sequences in the subcategories ${\rm mod}^{+\hspace{-3pt}}\La$ and ${\rm mod}^{-\hspace{-3pt}}\La$ of $\Mod\La$.

\begin{Lemma}\label{fdmod}

Let $\La=kQ/\mk I$ be locally semiprimary, where $Q$ is a quiver and $\mk I$ is a relation ideal in $kQ$. Then a module in $\Mod \La$ is finite dimensional if and only if it is finitely generated and finitely cogenerated.

\end{Lemma}

\noindent{\it Proof.} Let $M\in {\rm mod}^{\hspace{.4pt} b\hspace{-2.5pt}}\La.$ In particular, $M$ is finitely generated and $\soc M\in {\rm mod}^{\hspace{.4pt} b\hspace{-2.5pt}}\La.$ And there exist $x_1, \ldots, x_r\in Q_0$ such that $e_xM= 0$, for $x\in Q_0\backslash\{x_1, \ldots, x_r\}$. Let $0\ne m\in e_{x_i}M$ for some $1\le i\le r$. Suppose that $(\soc M)\cap (\La m)=0.$ By Proposition \ref{lfd_rad}(4), there exist $\alpha_n: a_n\to a_{n+1}$ in $Q_1$ with $a_1=x_i,$ such that $0\ne \bar \alpha_n \cdots \bar \alpha_1 m \in e_{a_n}M$, for all $n \ge 1.$ In particular, $a_n\in \{x_1, \ldots, x_r\}$ and $\alpha_n\cdots \alpha_1$ is a nonzero path in $Q$ from $x_i$ to $a_n$, for all $n\ge 1$. So, there exists some $1\le j\le r$ such that $Q$ has infinitely many nonzero paths from $x_i$ to $x_j,$ a contradiction. Thus, $\soc M$ is essential in $M$, and hence,  $M$ is finitely cogenerated; see (\ref{fd_soc}). 

Assume now that $M$ is finitely cogenerated, and finitely generated with an epimorphism $f: P_{x_1} \oplus \cdots \oplus P_{x_r} \to M,$ where $x_i\in Q_0$. Then \vspace{1pt} $M=\sum_{i=1}^r \La f(e_{x_i}),$ where $f(e_{x_i})\in e_{x_i}M,$ and by Lemma \ref{fd_soc}, $\soc M= e_{y_1}(\soc M) \oplus \cdots  \oplus e_{y_s}(\soc M)$, where $y_i\in Q_0$. For any $x\in Q_0$, ${\rm dim}(e_xM)\le \sum_{i=1}^r
{\rm dim}(e_x\La e_{x_i})<\infty$. Suppose that $e_xM\ne 0$ for some $x\in Q_0$. Then, $Q$ contains a path $\rho$ from some $x_i$ to $x$ such that $\bar{\rho} f(e_{x_i})\ne 0$. Since $\soc M$ is essential in $M$; see (\ref{fd_soc}), $Q$ contains a path $\eta$ from $x$ to some $y_j$ such that $\bar \eta (\bar \rho f(e_{x_i}))\ne 0.$ So, $\eta \rho$ is a nonzero path in $Q$ from $x_i$ to $y_j$ passing through $x$.
By Definition \ref{lspri_def}, $e_xM\ne 0$ for only finitely many $x\in Q_0$. Thus, $M\in {\rm mod}^{\hspace{.4pt} b\hspace{-2.5pt}}\La.$ The proof of the lemma is completed.

\begin{Remark}\label{fg_cfg_nfd} 

Lemma \ref{fdmod} does not hold for locally semiperfect algebras. Consider the locally semiperfect algebra $\La=kQ/\mk I$, where $Q$ is the quiver \vspace{-2pt} 
$$\xymatrixrowsep{20pt}\xymatrix{
 a_0\ar[r]^{\alpha_0} & a_1 \ar[r]^{\alpha_1} \ar[d]^{\gamma_1}  & a_2 \ar[d]^{\gamma_2}\ar[r]^{\alpha_2} & a_3 \ar[d]^{\gamma_3} \ar[r]^{\alpha_3} & a_4 \ar[d]^{\gamma_4} \ar@{.}[r] & \\
b_0 &b_1\ar[l]_{\beta_0} &b_2\ar[l]_{\beta_1} &b_3\ar[l]_{\beta_2} &b_4\ar[l]_{\beta_3} & \ar@{.}[l]} \vspace{-2pt}$$ and $\mk{I}$ is the ideal generated by $\beta_0\gamma_1\alpha_0-\beta_0\beta_1\cdots \beta_i \gamma_{i+1} \alpha_i \cdots \alpha_1\alpha_0$ with $i\ge 1,$ and $\beta_i\gamma_{i+1} \alpha_i-\beta_i\beta_{i+1}\cdots \beta_j\gamma_{j+1} \alpha_j \cdots \alpha_{i+1}\alpha_i$ with $j>i\ge 1$. Note that $P_{a_0}\cong I_{b_0}$. Thus, $P_{a_0}$ is finitely generated, finitely cogenerated and infinite dimensional.

\end{Remark}

The following result is handy for constructing examples of locally semiprimary algebras, which are also locally left or right noetherian.

\begin{Lemma} 

Let $\La=kQ/\mk{I}$, where $Q$ is a quiver and $\mk{I}$ is a relation ideal in $kQ$. If $\La$ is locally left or right eventually multiserial, then $\La$ is locally semiprimary if and only if all oriented cycles in $Q$ are nilpotent in $\La$.

\end{Lemma}

\noindent{\it Proof.} We consider only the case where $\La$ is locally left eventually multiserial. The necessity is evident. Assume that all oriented cycles in $Q$ are nilpotent in $\La$. Suppose that $Q$ has infinitely many nonzero paths from $x$ to $y$, for some $x, y\in Q_0$. Since $Q$ is left locally finite, it contains an infinite path \vspace{-3pt} $$\zeta: \; \xymatrixcolsep{20pt}\xymatrix{x=y_0\ar[r]^-{\alpha_1} & y_1\ar[r]^-{\alpha_2} & y_2 \ar[r]^{\alpha_3} & \cdots \ar[r]^{\alpha_i} & y_i \ar[r]^{\alpha_{i+1}} & \cdots } \vspace{-3pt} $$ such that $Q$ has a path $\rho_i: y_i\rightsquigarrow y$ such that $\rho_i \alpha_i\cdots \alpha_1$ is a nonzero path, for all $i\ge 1$. Let $n_x$ be an integer as in Definition \ref{mserial}(1). Then, for any $i\ge n_x$, the infinite subpath $\zeta_i$ of $\zeta$ starting with $y_i$ starts with $\rho_i$. So, we may assume that $y_{n_x}=y$. It is not difficult to see that $\zeta_{\hspace{.5pt} n_x}$ starts with $\rho_{n_x}^s$ for infinitely many $s>0$, contrary to the assumption. The proof of the lemma is completed.

\subsection{\sc Finite presented and finitely copresented modules} We begin with the following statement, which is a special case of Proposition \ref{C-proj}.

\begin{Prop}\label{lfd_fpres}

Let $\La=kQ/\mk I$ be locally semiperfect, where $Q$ is a quiver and $\mk I$ is a relation ideal in $kQ$. Then

\begin{enumerate}[$(1)$]

\vspace{-2pt} 

\item ${\rm mod}^{+\hspace{-3pt}}\La$ and ${\rm mod}^{-\hspace{-3pt}}\La$ are Hom-finite Krull-Schmidt subcategories of $\mmod \La\hspace{.5pt};$ 

\item The duality $\mk D: \mmod \La^\circ \to \mmod \La$ restricts to a duality $\mk D: {\rm mod}^{+\hspace{-3pt}}\La^\circ \to {\rm mod}^{-\hspace{-3pt}}\La$.

\end{enumerate} \end{Prop}


We shall need the property that finite dimensional modules are all finitely presented or all finitely copresented.

\begin{Lemma}\label{fd_fp}

Let $\La=kQ/\mk I$ be locally semiperfect, where $Q$ is a quiver and $\mk I$ is a relation ideal in $kQ$. Then ${\rm mod}^{\hspace{.4pt}b\hspace{-2.5pt}}\La \subseteq {\rm mod}^{+\hspace{-3pt}}\La$ if and only if $Q$ is left locally finite$\,;$ and ${\rm mod}^{\hspace{.4pt}b\hspace{-2.5pt}}\La \subseteq {\rm mod}^{-\hspace{-3pt}}\La$ if and only if $Q$ is right locally finite.  \end{Lemma}

\vspace{-1.5pt} 

\noindent{\it Proof.} Since ${\rm mod}^{+\hspace{-3pt}}\La$ is extension-closed in $\ModLa$, by Proposition \ref{lfd_rad}(3), we see that ${\rm mod}^{\hspace{.4pt}b\hspace{-2.5pt}}\La \subseteq {\rm mod}^{+\hspace{-3pt}}\La$ if and only if  
$S_x\in {\rm mod}^{+\hspace{-3pt}}\La$ for all $x\in Q_0$. And for each $x\in Q_0$, there exists a short exact sequence \vspace{-.5pt} 
$\hspace{-2pt}\xymatrixcolsep{18pt}\xymatrix{0\ar[r] & Je_x\ar[r] & P_x\ar[r] & S_x\ar[r] & 0}\hspace{-2pt}$ in $\Mod\La$. Since $\{\bar \alpha \mid \alpha \in Q_1(x, -)\}$ is a minimal generating set for $Je_x$, by Schanul's Lemma, $S_x\in {\rm mod}^{+\hspace{-3pt}}\La$ if and only if $\{\bar \alpha \mid \alpha \in Q_1(x, -)\}$ is finite, that is, $Q_1(x, -)$ is finite. 
The first part of the lemma holds. Next, by Propositions \ref{duality_lfdm} and \ref{lfd_fpres}(2), ${\rm mod}^{\hspace{.4pt}b\hspace{-2.5pt}}\La \subseteq {\rm mod}^{-\hspace{-3pt}}\La$ if and only if ${\rm mod}^{\hspace{.4pt}b\hspace{-2.5pt}}\La^\circ \subseteq {\rm mod}^{+\hspace{-3pt}}\La^\circ$, or equivalently, $Q^\circ$ is left locally finite, that is, $Q$ is right locally finite. The proof of the lemma is completed. 

\vspace{2pt}

We shall also need the following easy statement.

\begin{Lemma}\label{llb_lrb}

Let $\La=kQ/R$ be semiperfect, where $Q$ is a quiver and $\mk I$ is a relation ideal in $kQ$. Then, ${\rm mod}^{+\hspace{-3pt}}\La={\rm mod}^{\hspace{.4pt}b\hspace{-3pt}}\La$ if and only if $\La$ is locally left bounded$\hspace{.5pt};$ and ${\rm mod}^{-\hspace{-3pt}}\La={\rm mod}^{\hspace{.4pt}b\hspace{-3pt}}\La$ if and only if $\La$ is locally right bounded.

\end{Lemma}

\noindent{\it Proof.} We shall only prove the sufficiency of the first part of the lemma. Assume that $\La$ is locally left bounded. Then, ${\rm proj}\La \subseteq {\rm mod}^{\hspace{.4pt}b\hspace{-3pt}}\La,$ and hence, ${\rm mod}^{+\hspace{-3pt}}\La\subseteq {\rm mod}^{\hspace{.4pt}b\hspace{-3pt}}\La$. Moreover, $Q$ is clearly left locally finite, and hence, ${\rm mod}^{\hspace{.4pt}b\hspace{-3pt}}\La\subseteq {\rm mod}^{+\hspace{-3pt}}\La;$ see (\ref{fd_fp}). The proof of the lemma is completed.

\vspace{2pt}

Clearly, the projective modules in ${\rm mod}^{+\hspace{-3pt}}\La$ are ext-projective, and the injective modules in ${\rm mod}^{-\hspace{-3pt}}\La$ are ext-injective. The converses are not true in general. 

\begin{Lemma}\label{ext_proj_ab}

Let $\La=kQ/R,$ where $Q$ is a quiver and $\mk I$ is a relation ideal in $kQ$. 

\begin{enumerate}[$(1)$]

\vspace{-4pt}

\item If ${\rm mod}^{+\hspace{-3pt}}\La$ is abelian, then its ext-projective objects are the modules in ${\rm proj}\La.$

\vspace{1pt}

\item If ${\rm mod}^{-\hspace{-3pt}}\La$ is abelian, then its ext-injective objects are the modules in ${\rm inj}\La.$

\end{enumerate} \end{Lemma}

\noindent{\it Proof.} We shall only prove Statement (1). Suppose that ${\rm mod}^{+\hspace{-3pt}}\La$ be abelian. Consider a nonprojective module $M$ in ${\rm mod}^{+\hspace{-3pt}}\La.$ \vspace{-.5pt} Then there exists a nonsplit short exact sequence $\xymatrixcolsep{20pt}\xymatrix{\hspace{-2pt} 0\ar[r] & L\ar[r] & P\ar[r] & M\ar[r] & 0}\hspace{-2pt}\vspace{-1pt}$ in $\Mod\La$ where $P\in {\rm proj}\La$ and $L$ is finitely generated. By Proposition \ref{fpres_AB}, $L\in {\rm mod}^{+\hspace{-3pt}}\La,$ and hence, $M$ is not ext-projective in ${\rm mod}^{+\hspace{-3pt}}\La$. The proof of the lemma is completed.

\vspace{2pt}

We provide some classes of algebras $\La$ for which ${\rm mod}^{+\hspace{-3pt}}\La$ or ${\rm mod}^{-\hspace{-3pt}}\La$ is abelian.

\begin{Prop} \label{fpres_ab} Let $\La=kQ/\mk I$ be locally semiperfect, where $Q$ is a quiver and $\mk I$ is a relation ideal in $kQ$.

\begin{enumerate}[$(1)$] 

\vspace{-3.5pt}

\item If $\mk I=0$, then ${\rm mod}^{+\hspace{-3pt}}\La$ and ${\rm mod}^{-\hspace{-3pt}}\La$ are abelian.

\vspace{.5pt}

\item If $\La$ is locally left $(\hspace{-.5pt}$resp. right$\hspace{.8pt})$ noetherian, then ${\rm mod}^{+\hspace{-3pt}}\La \hspace{.4pt}(\hspace{-.5pt}$resp. ${\rm mod}^{-\hspace{-3pt}}\La)$ is abelian.

\end{enumerate} \end{Prop}

\noindent {\it Proof.} (1) Assume that $\mk I=0$. Then $\La$ is hereditary; see \cite[(8.2)]{GaR}. Thus,  ${\rm proj}\La$ is closed under finitely generated submodules and ${\rm inj}\La$ is closed under finitely cogenerated quotient-modules. By Proposition \ref{fpres_AB}, ${\rm mod}^{+\hspace{-3pt}}\La$ and ${\rm mod}^{-\hspace{-3pt}}\La$ are abelian.

(2) Suppose that $\La$ is locally left noetherian. Then the modules in ${\rm proj}\La$ are noetherian. Thus, finitely generated modules are finitely presented. By Proposition \ref{fpres_AB}, ${\rm mod}^{+\hspace{-3pt}}\La$ is abelian. Assume that $\La$ is locally right noetherian. Then $\La^\circ$ is locally left noetherian, and hence, ${\rm mod}^{+\hspace{-3pt}}\La^\circ$ is abelian. By Proposition \ref{lfd_fpres}(2), ${\rm mod}^{-\hspace{-2.5pt}}\La$ is abelian. The proof of the proposition is completed.

\begin{Exam}\label{Ab_subcat_ex} Consider $\La=kQ/\mk I$, where $Q$ is the quiver \vspace{-6pt} 
$$\hspace{-20pt}\xymatrixrowsep{20pt}\xymatrix{
 \ar@{.}[r] & a_3 \ar[r]^{\alpha_3} & a_2\ar[r]^{\alpha_2} & a_1 \ar[r]^{\alpha_{1}}  & a_0\ar@<.5ex>[r]^\alpha\ar@<-.6ex>[r]_\beta& a_{-1} \ar[r]\ar[d] & a_{-2} \ar[d] \ar[r]& a_{-3} \ar@<-2ex>[d]\ar@{.}[r] & \\
& b_3\ar[u]_{\beta_3} &  b_2\ar[u]_{\beta_2}&b_1\ar[u]_{\beta_1}&&c_1&c_1& \hspace{-15pt} c_2 } \vspace{-3pt}$$
and $\mk{I}=\langle \alpha\alpha_1, \beta\alpha_1, \alpha_i\beta_i \mid i\ge 1\rangle.$ Then, $\La$ is locally right eventually multiserial. By Propositions \ref{mul-noe} and \ref{fpres_ab}, 
${\rm mod}^{-\hspace{-3pt}}\La$ is abelian. Although $\La$ is neither locally left noetherian nor hereditary, we claim that
${\rm mod}^{+\hspace{-3pt}}\La$ is abelian. 

Indeed let $\La'=kQ'$, where $Q'$ is the full subquiver of $Q$ generated by the $a_i$ and $c_j$ with $i\le 0$ and $j\ge 1$. Then ${\rm proj} \La'\subseteq {\rm proj}\La$. Let $M$ be a finitely generated submodule of some module $P$ in ${\rm proj}\La$. Write $P'$ for the maximal direct summand of $P$ lying in ${\rm proj}\La'$. Note that $M'=\sum_{a_0\ne x\in Q'_0} e_xM$ is a finitely generated submodule of $P'\cap M$ such that $M/M'\in {\rm mod}^{\hspace{.5pt}b\hspace{-2.5pt}}\La \subseteq {\rm mod}^+\hspace{-3pt}\La;$ see (\ref{fd_fp}). 
Since $\La'$ is hereditary, $M'\in {\rm proj}\La'\subseteq {\rm mod}^+\hspace{-3pt}\La$. Therefore, $M\in {\rm mod}^+\hspace{-3pt}\La$. Now it follows from Proposition \ref{fpres_AB}(1) that ${\rm mod}^+\hspace{-3pt}\La$ is abelian.

\end{Exam}

\vspace{-6pt}

\subsection{\sc Almost split sequences} We start with the locally semiperfect case. Since ${\rm mod}^{+\hspace{-3pt}}\La$ and ${\rm mod}^{-\hspace{-3pt}}\La$ are Krull-Schmidt; see (\ref{lfd_fpres}), as a special case of Theorem \ref{C_mod_ass}, we immediately obtain a general existence theorem for $\Mod\La$ as follows.

\begin{Theo}\label{ass_lfd}

Let $\La=kQ/\mk{I}$ be a locally semiperfect algebra, where $Q$ is a quiver and $\mk I$ is a relation ideal in $kQ$. Consider an indecomposable module $M\in \Mod\La.$

\begin{enumerate}[$(1)$]

\vspace{-2pt}

\item If $M\in \mmod^{+\hspace{-3pt}}\La$ \vspace{-.5pt} is not projective, then there exists an almost split sequence $\xymatrixcolsep{20pt}\xymatrix{ 0\ar[r] &\tau M\ar[r] & N \ar[r] &M \ar[r] &0}\hspace{-2pt}$ in $\ModLa,$ where $\tau M\in \mmod^{-\hspace{-3pt}}\La.$

\item If $M \in \mmod^{-\hspace{-3pt}}\La$ is not injective, \vspace{-1.2pt} then there exists an almost split sequence $\xymatrixcolsep{20pt}\xymatrix{0\ar[r] & M\ar[r] & N \ar[r] &\tau^-\hspace{-1.2pt}M \ar[r] &0}\hspace{-2pt}$ in $\ModLa,$ where $\tau^-\hspace{-1.2pt} M\in \mmod^{+\hspace{-3pt}}\La.$

\end{enumerate} \end{Theo}

\begin{Remark} 
(1) In case $Q$ is strongly locally finite and $\mk I=0,$ Theorem \ref{ass_lfd} strengthens the result of Theorem 2.8 in \cite{BLP}.

\vspace{.8pt}

\noindent (2) In case $\La$ is locally left or right bounded, Theorem \ref{ass_lfd} strengthens the result of Theorem 3.4.1 in \cite{Gup}. Note that the latter is obtained using the classical approach, and contains an error claiming that $\tau M$ and $\tau^-M$ are finite dimensional. 

\end{Remark}

Next, we shall concentrate on almost split sequences in ${\rm mod}^{+\hspace{-3pt}}\La$ and ${\rm mod}^{-\hspace{-3pt}}\La$ in the locally semiprimary case. This will be based on the following description of irreducible morphisms in ${\rm mod}^{+\hspace{-3pt}}\La$ and ${\rm mod}^{-\hspace{-3pt}}\La,$ which covers the result in \cite[(3.1)]{BLP}. 

\begin{Lemma}\label{ker_irr_epi}

Let $\La=kQ/\mk I$ be locally semiprimary, where $Q$ is a quiver and $\mk I$ is a relation ideal in $kQ$. 

\begin{enumerate}[$(1)$]

\vspace{-1pt}

\item If $Q$ is left locally finite, then every irreducible epimorphism in ${\rm mod}^{+\hspace{-3pt}}\La$ has a finite dimensional kernel.

\vspace{1pt}

\item If $Q$ is right locally finite, then every irreducible monomorphism in ${\rm mod}^{-\hspace{-3pt}}\La$ has a finite dimensional cokernel.

\end{enumerate}

\end{Lemma}

\vspace{-1.5pt}

\noindent{\it Proof.} We shall only prove Statement (1). Assume that $Q$ is left locally finite. By Lemma \ref{fd_fp}, ${\rm mod}^{\hspace{.4pt}b}\hspace{-3pt}\La \subseteq {\rm mod}^{+}\hspace{-3pt}\La.$ Consider an irreducible epimorphism $f: M\to N$ in ${\rm mod}^{+\hspace{-3pt}}\La$. We may assume that $N$ is indecomposable. In view of Theorem \ref{ass_lfd}(1), we can construct a commutative diagram with exact rows \vspace{-7pt}
$$\xymatrixrowsep{22pt}\xymatrix{0\ar[r]& L \ar[r]^g \ar[d]^u & M \ar[d]^v\ar[r]^f & N\ar[r] \ar[r] \ar@{=}[d] & 0\\ 0\ar[r]& X \ar[r]^h & Y \ar[r]^p& N \ar[r] & 0,} \vspace{-8pt}$$ where the bottom row is an almost split sequence in $\ModLa$ with $X\in {\rm mod}^{-\hspace{-3pt}}\La$. Writing $X'={\rm Im}(u)$ and $Y'={\rm Im}(v),$ we obtain factorizations $u=ju'$ and $v=qv',$ where $u': L\to X'$
and $v': M\to Y'$ are epimorphisms; $j: X'\to X$ and $q: Y'\to Y$ are inclusion maps. This yields a commutative diagram with exact rows 
\vspace{-1pt}
$$\xymatrixrowsep{20pt}\xymatrix{0\ar[r]& L \ar[r]^g \ar[d]^{u'} & M \ar[d]^{v'}\ar[r]^f & N\ar[r] \ar[r] \ar@{=}[d] & 0\\ 0\ar[r]& X' \ar[r]^{h'} & Y' \ar[r]^{pq}& N \ar[r] & 0.} \vspace{-0pt}$$ 

Since $L$ is finitely generated, so is $X'$. And since $X$ is finitely cogenerated, so is $X';$ see (\ref{fd_soc}). Thus, $X'\in {\rm mod}^{\hspace{.4pt}b}\hspace{-3pt}\La \subseteq {\rm mod}^{+}\hspace{-3pt}\La;$ see (\ref{fdmod}), and hence, $Y'\in {\rm mod}^{+}\hspace{-3pt}\La;$ see (\ref{subcat_pres}). Since $p$ is not a retraction, neither is $pq: Y'\to N$. Thus, $v'$ is a section. Then, $v'$ is an isomorphism, and so is $u'$. So, $L$ is finite dimensional. The proof of the lemma is completed.

\vspace{2pt}

We are now ready to describe the almost split sequences in ${\rm mod}^{+\hspace{-3pt}}\La$ which, in the hereditary case, strengthens the result stated in \cite[(3.6)]{BLP}. 

\begin{Theo}\label{ass_fpres}

Let $\La=kQ/\mk I$ be a locally semiprimary algebra, where $Q$ is a locally finite quiver and $\mk I$ is a relation ideal in $kQ$. \vspace{.5pt} Then the almost split sequences in ${\rm mod}^{+\hspace{-3pt}}\La$ are those in $\ModLa$ with a finite dimensional starting term. Moreover,

\begin{enumerate}[$(1)$]

\vspace{-.5pt}

\item an indecomosable nonprojective module $M$ in ${\rm mod}^{+\hspace{-3pt}}\La$ is the ending term of an almost split sequence in ${\rm mod}^{+\hspace{-3pt}}\La$ if and only if $\tau M$ is finite dimensional$\hspace{.5pt};$

\vspace{1pt}

\item an indecomosable noninjective module $N$ in ${\rm mod}^{+\hspace{-3pt}}\La$ is the starting term of an almost split sequence in ${\rm mod}^{+\hspace{-3pt}}\La$ if and only if $N$ is finite dimensional.

\end{enumerate}

\end{Theo}

\noindent{\it Proof.} Since $Q$ is locally finite, ${\rm mod}^{\hspace{.5pt}b\hspace{-2.8pt}}\La \hspace{-.4pt}= \hspace{-.4pt} {\rm mod}^{+\hspace{-3pt}}\La \hspace{.6pt} \cap \hspace{.6pt} {\rm mod}^{-\hspace{-3pt}}\La;$ 
see (\ref{fdmod}) and (\ref{fd_fp}). 
Consider an almost split sequence $\xymatrixcolsep{20pt}\xymatrix{0\ar[r]& X \ar[r] & Y \ar[r] & Z \ar[r] & 0}\hspace{-2pt}$ in $\ModLa,$ where $X$ is finite dimensional. By Theorem \ref{ass_lfd}(2), $Z\cong \tau^-X\in {\rm mod}^{+\hspace{-3pt}}\La.$ Since ${\rm mod}^{+\hspace{-3pt}}\La$ is extension-closed in $\ModLa$, this is an almost split sequence in ${\rm mod}^{+\hspace{-3pt}}\La.$ Conversely, let $\, (*) \hspace{5pt} \xymatrixcolsep{20pt}\xymatrix{0\ar[r]& M \ar[r] & N \ar[r] & L \ar[r] & 0}\hspace{-2pt}$ be an almost split sequence in ${\rm mod}^{+\hspace{-3pt}}\La.$ By Lemma \ref{ker_irr_epi}(1), $M\in {\rm mod}^{\hspace{.4pt} b\hspace{-2.5pt}}\La,$  \vspace{-1.5pt} and by Theorem \ref{ass_lfd}(2), there exists an almost split sequence $\xymatrixcolsep{20pt}\xymatrix{\hspace{-2pt} 0\ar[r]& M \ar[r] & E \ar[r] & \tau^-\hspace{-1pt}M \ar[r] & 0}\hspace{-2pt}\vspace{-.5pt}$ in $\ModLa,$ where $\tau^-\hspace{-1pt}M\in {\rm mod}^{+\hspace{-3pt}}\La.$ Since $M\in {\rm mod}^{+\hspace{-3pt}}\La,$ this is an almost split sequence in ${\rm mod}^{+\hspace{-3pt}}\La$, and hence, it is isomorphic to $(*)$. 
This establishes the first part of the theorem. Combining this part with Theorem \ref{ass_lfd}, we deduce easily Statements (1) and (2). The proof of the theorem is completed.

\begin{Remark} In view of Theorem \ref{ass_fpres}(2), we see that the finite dimensional ext-injective objects in ${\rm mod}^{+\hspace{-3pt}}\La$ are the finite dimensional modules in ${\rm inj}\La.$ \end{Remark}

Dually, we can describe the almost split sequences in  ${\rm mod}^{-\hspace{-3pt}}\La$ as follows.

\begin{Theo}\label{ass_fcp}

Let $\La=kQ/\mk I$ be a locally semiprimary algebra, where $Q$ is a locally finite quiver and $\mk I$ is a relation ideal in $kQ$. Then the almost split sequences in ${\rm mod}^{-\hspace{-3pt}}\La$ are those in $\ModLa$ with a finite dimensional ending term. Moreover, 

\begin{enumerate}[$(1)$]

\vspace{-1pt}

\item an indecomposable noninjective module $N$ in ${\rm mod}^{-\hspace{-3pt}}\La$ is the starting term of an almost split sequence in ${\rm mod}^{-\hspace{-3pt}}\La$ if and only if $\tau^-N$ is finite dimensional$\,;$

\vspace{.5pt}

\item an indecomposable nonprojective module $M$ in ${\rm mod}^{-\hspace{-3pt}}\La$ is the ending term of an almost split sequence in ${\rm mod}^{-\hspace{-3pt}}\La$ if and only if  $M$ is finite dimensional.

\end{enumerate}\end{Theo}

\begin{Remark}
 
In view of \ref{ass_fcp}(2), we see that the finite dimensional ext-projective objects in ${\rm mod}^{-\hspace{-3pt}}\La$ are the finite dimensional modules in ${\rm proj}\La.$ 

\end{Remark}


\subsection{\sc Subcategories having almost split sequences} The objective of this subsection is to study when ${\rm mod}^{+\hspace{-3pt}}\La$, ${\rm mod}^{-\hspace{-3pt}}\La$ and ${\rm mod}^{\hspace{.5pt}b\hspace{-2.8pt}}\La$ have almost split sequences on one or both sides in case $\La$ is locally semi\-primary given by a locally finite quiver.


\begin{Theo}\label{ass_fp_cat}

Let $\La=kQ/R$ be a locally semiprimary algebra, where $Q$ is a locally finite quiver and $\mk I$ is a relation ideal in $kQ$. 

\begin{enumerate}[$(1)$]

\vspace{-1pt}

\item If  $\La$ is locally left bounded, then ${\rm mod}^{+\hspace{-3pt}}\La$ has almost split sequences on the left. And the converse holds in case ${\rm mod}^{+\hspace{-3pt}}\La$ is abelian.

\vspace{.5pt}

\item If the indecomposable noninjective modules in ${\rm mod}^{-\hspace{-3pt}}\La$ are finite dimensional $($in particular, if $\La$ is locally right bounded$\hspace{.8pt})$, then ${\rm mod}^{+\hspace{-3pt}}\La$ has almost split sequences on the right. And the converse holds in case ${\rm mod}^{+\hspace{-3pt}}\La$ is abelian.

\end{enumerate}\end{Theo}

\vspace{-1pt}

\noindent{\it Proof.} (1) Suppose that $\La$ is locally left bounded. Then,   ${\rm mod}^{+\hspace{-3pt}}\La={\rm mod}^{\hspace{.4pt}b\hspace{-3pt}}\La$ by Lemma \ref{llb_lrb}. If $M\in {\rm mod}^{+\hspace{-3pt}}\La$ is indecomposable and not ext-injective, then it is not injective, and by Theorem \ref{ass_fpres}(2), $M$ is the starting term of an almost split sequence in ${\rm mod}^{+\hspace{-3pt}}\La.$ Thus, ${\rm mod}^{+\hspace{-3pt}}\La$ has almost split sequences on the left.

Suppose now that  $\La$ is not locally left bounded and ${\rm mod}^{+\hspace{-3pt}}\La$ is abelian. Then $P_x$ is infinite dimensional for some $x\in Q_0$. Since $Q$ is locally finite, 
by Proposition \ref{fpres_AB}(1), $\rad P_x\in {\rm mod}^{+\hspace{-3pt}}\La.$ Let $N$ be an infinite dimensional indecomposable direct summand of $\rad P_x$. \vspace{-1pt} Then, we obtain a nonsplit short exact sequence 
$\xymatrixcolsep{20pt}\xymatrix{\hspace{-3pt} 0\ar[r] & N\ar[r] &  P_x\ar[r] & P_x/N \ar[r] & 0}\vspace{-2pt}$ in ${\rm mod}^{+\hspace{-3pt}}\La$. In particular, $N$ is not ext-injective in ${\rm mod}^{+\hspace{-3pt}}\La$. By Theorem \ref{ass_fpres}(2), $N$ is not the starting term of any almost split sequence in ${\rm mod}^{+\hspace{-3pt}}\La$. Thus, ${\rm mod}^{+\hspace{-3pt}}\La$ does not have almost split sequences on the left. 

(2) Suppose that all indecomposable noninjective modules in ${\rm mod}^{-\hspace{-3pt}}\La$ are finite dimensional. If $M\in {\rm mod}^{+\hspace{-3pt}}\La$ is indecomposable and not ext-projective, then $M$ is not projective, and by Theorem \ref{ass_lfd}(1), \vspace{-1pt} there exists an almost split sequence $\xymatrixcolsep{20pt}\xymatrix{\hspace{-3pt}0\ar[r] & \tau M \ar[r] & N \ar[r] & M\ar[r] & 0\hspace{-1pt}}\vspace{-1pt}$ in $\Mod\La,$ where $\tau M\in {\rm mod}^{-\hspace{-3pt}}\La.$ Since $\tau M\in {\rm mod}^{\hspace{.5pt}b\hspace{-2.8pt}}\La$ by the assumption, this is an almost split sequence in ${\rm mod}^{+\hspace{-3pt}}\La$. Therefore, ${\rm mod}^{+\hspace{-3pt}}\La$ has almost split sequences on the right. 

\vspace{.5pt}

Assume now that ${\rm mod}^{+\hspace{-3pt}}\La$ is abelian and ${\rm mod}^{-\hspace{-3pt}}\La$ contains an infinite dimensional indecomposable noninjective module $M$. \vspace{-1pt} By Theorem \ref{ass_lfd}(2), there exists an almost split sequence $\xymatrixcolsep{20pt}\xymatrix{0\ar[r] & M \ar[r] & N  \ar[r] & L\ar[r] & 0\hspace{-2pt}}$ \vspace{-.5pt} in  $\Mod\La,$ where $L\in {\rm mod}^{+\hspace{-3pt}}\La$. Since $\tau L \cong M,$ by Theorem \ref{ass_fpres}(1),  ${\rm mod}^{+\hspace{-3pt}}\La$ has no almost split sequence ending with $L$. Since $L$ is not ext-projective in ${\rm mod}^{+\hspace{-3pt}}\La$ by Lemma \ref{ext_proj_ab}(1), ${\rm mod}^{+\hspace{-3pt}}\La$ does not have almost split sequences on the right. The proof of the theorem is completed.

\vspace{-1pt}

\begin{Remark}

In case $\mk I=0$, we see from Lemma \ref{fpres_ab}(1) that Theorem \ref{ass_fp_cat} covers the result in \cite[(3.7)]{BLP}; compare also \cite[(3.6.1)]{LL25}.

\end{Remark}

Dually, we have the following statement.

\vspace{-2pt}

\begin{Theo}\label{ass_fcp_cat}

Let $\La=kQ/R$ be a locally semiprimary algebra, where $Q$ is a locally finite quiver and $\mk I$ is a relation ideal in $kQ$. 

\begin{enumerate}[$(1)$]

\vspace{-1.5pt}

\item If $\La$ is locally right bounded, then ${\rm mod}^{-\hspace{-3pt}}\La$ has almost split sequences on the right. And the converse holds if ${\rm mod}^{-\hspace{-3pt}}\La$ is abelian.

\vspace{.5pt}

\item If the indecomposable nonprojective modules in ${\rm mod}^{+\hspace{-3pt}}\La$ are finite dimensional $($in particular, if $\La$ is locally left bounded$\hspace{.5pt})$, then ${\rm mod}^{-\hspace{-3pt}}\La$ has almost split sequences on the left. And the converse holds if ${\rm mod}^{-\hspace{-3pt}}\La$ is abelian.
    
\end{enumerate}\end{Theo}

\begin{Exam}\label{ex_ass_fp_cat} (1) Consider $\La=kQ/\mk I$, where $Q$ is the quiver 
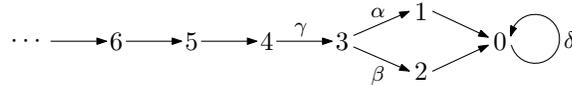
\vspace{-2pt}
$$\begin{tikzpicture}[-{Stealth[inset=0pt,length=4.5pt,angle'=35,round,bend]},scale=.5]
\node at (15.95, 0.05){\small $\delta$};
\node at (14.1,0) {$0$};
\draw (14.4, -0.1) arc (14.5:355:-0.65);
\draw  (12.3, 0.8) -- (13.8,0.1);
\draw  (12.3, -0.8) -- (13.8,-0.1);
\node at (10.85, 0.75) {\footnotesize $\alpha$};
\draw  (10.2, 0.15) -- (11.7,0.8);
\node at (12,0.8) {$1$};
\draw  (10.2,-0.15) -- (11.7,-0.8);
\node at (12,-0.8) {$2$};
\node at (10.85,-0.9){\footnotesize $\beta$};
\node at (9.9,0){$3$};
\draw  (8.15,0) -- (9.65,0);
\node at (8.8,0.3){\footnotesize $\gamma$};
\node at (7.9,0){$4$};
\draw  (6.15,0) -- (7.65,0);
\node at (5.9,0) {$5$};
\draw (4.15,0) -- (5.65,0);
\node at (3.9,0) {$6$};
\draw (2.15,0) -- (3.65,0);
\node at (1.55,0){$\cdots$};  
\end{tikzpicture}\vspace{-8.5pt} $$
and $\mk I=\langle \alpha \gamma, \beta \gamma, \delta^3\hspace{.4pt}\rangle$. Then, $\La$ is locally eventually multiserial. By Lemmas \ref{mul-noe} and \ref{fpres_ab},  ${\rm mod}^{-\hspace{-3pt}}\La$ and ${\rm mod}^{+\hspace{-3pt}}\La$ are abelian. Since $\La$ is locally left bounded but not locally right bounded, by Theorem \ref{ass_fcp_cat}, ${\rm mod}^{-\hspace{-3pt}}\La$ has almost split sequences on the left but not on the right. And by Theorem \ref{ass_fp_cat}(1), ${\rm mod}^{+\hspace{-3pt}}\La$ has almost split sequences on the left. We claim that it also has almost split sequences on the right. 

Indeed let $\La'=kQ'$ and $\La''=(kQ'')/\langle \delta^3\hspace{.4pt}\rangle$, where $Q'$ and $Q''$ are the full subquivers of $Q$ generated by the vertices $i\ge 3$ and by the vertices $j$ with $0\le j\le 3$, respectively. Then, every module $M$ in $\ModLa$ is decomposed as $M=M'\oplus M''$, where $M'\in \ModLa'$ and $M''\in \ModLa''.$ Thus, if $M$ is an infinite dimensional indecomposable module in ${\rm mod}^{-\hspace{-3pt}}\La$, then $M\in {\rm mod}^{-\hspace{-3pt}}\La',$ and hence, $M\cong I_i$ for some $i\ge 3$. By Theorem \ref{ass_fp_cat}(2), ${\rm mod}^{+\hspace{-3pt}}\La$ has almost split sequences on the right.

\vspace{4pt}

\noindent (2) Consider $\La=kQ/\mk I$, where $Q$ is the quiver \vspace{-4pt} 
$$\hspace{0pt}\xymatrixrowsep{20pt}\xymatrix{
 \ar@{.}[r] & a_3 \ar[r]^{\alpha_3} & a_2\ar[r]^{\alpha_2} & a_1 \ar[r]^{\alpha_{1}}  & a_0\ar@<.5ex>[r]^\alpha\ar@<-.6ex>[r]_\beta& a_{-1} \ar[r]\ar[d] & a_{-2} \ar[d] \ar[r]& a_{-3} \ar[d] \ar@{.}[r] & \\
&b_3\ar[u]_{\beta_3} &  b_2\ar[u]_{\beta_2}&b_1\ar[u]_{\beta_1}&&c_1&c_1& c_2 } \vspace{-3pt}$$
and $\mk{I}=\langle \alpha\alpha_1, \beta\alpha_1, \alpha_i\beta_i \mid i\ge 1\rangle.$ We claim that 
neither ${\rm mod}^{+\hspace{-3pt}}\La$ nor ${\rm mod}^{-\hspace{-3pt}}\La$ 
has almost split sequences on either side. As seen in Example \ref{Ab_subcat_ex}, ${\rm mod}^{+\hspace{-3pt}}\La$ and ${\rm mod}^{-\hspace{-3pt}}\La$ are abelian. 
Since $\La$ is not left locally bounded, by Theorem \ref{ass_fp_cat}(1),
${\rm mod}^{+\hspace{-3pt}}\La$ does not have almost split sequences on the left. On the other hand, let $M$ be the submodule of $I_{a_1}$ such that $I_{a_1}/M\cong S_{b_1}.$ Then $M$ is an infinite dimensional indecomposable noninjective module in ${\rm mod}^{-\hspace{-3pt}}\La$. By Theorem \ref{ass_fp_cat}(2),
${\rm mod}^{+\hspace{-3pt}}\La$ does not have almost split sequences on the right. Similarly, we see from Theorem \ref{ass_fcp_cat} that ${\rm mod}^{-\hspace{-3pt}}\La$ does not have almost split sequences on either side.

\end{Exam}

\vspace{-1pt}

To conclude the paper, we study when ${\rm mod}^{\hspace{.5pt}b\hspace{-2.5pt}}\La$ has almost split sequences. 

\vspace{-1pt}

\begin{Theo}\label{ass_fdm_cat}

Let $\La=kQ/R$ be a locally semiprimary algebra, where $Q$ is a locally finite quiver and $\mk I$ is a relation ideal in $kQ$. 

\begin{enumerate}[$(1)$]

\vspace{-1pt}
    
\item If the indecomposable noninjective modules in ${\rm mod}^{-\hspace{-3pt}}\La$ are finite dimensional, then ${\rm mod}^{\hspace{.5pt}b\hspace{-2.5pt}}\La$ has almost split sequences on the right.

\item If the indecomposable nonprojective modules in ${\rm mod}^{+\hspace{-3pt}}\La$ are finite dimensional, then ${\rm mod}^{\hspace{.5pt}b\hspace{-2.5pt}}\La$ has almost split sequences on the left.

\item In each of the above two situations, the almost split sequences in ${\rm mod}^{\hspace{.5pt}b\hspace{-2.5pt}}\La$ are almost split sequences in $\ModLa$.

\end{enumerate}\end{Theo}

\noindent{\it Proof.} We shall only consider the case where all indecomposable noninjective mo\-dules in ${\rm mod}^{-\hspace{-3pt}}\La$ are finite dimensional. Let $M\in {\rm mod}^{\hspace{.5pt}b\hspace{-2.5pt}}\La$ be indecomposable and not ext-projective. By Lemma \ref{fd_fp}, \vspace{-.5pt} $M\in {\rm mod}^{+\hspace{-3pt}}\La$ is not projective, and by Theorem \ref{ass_lfd}, $\ModLa$ has an almost split sequence $\hspace{-1pt}\xymatrixcolsep{20pt}\xymatrix{0\ar[r] & \tau M \ar[r] & N\ar[r] & M\ar[r] & 0}\hspace{-1pt}\vspace{-1pt}$ with $\tau M\in {\rm mod}^{-\hspace{-3pt}}\La.$ Since $\tau M\in {\rm mod}^{\hspace{.5pt}b\hspace{-2.5pt}}\La,$  this is an almost split sequence in ${\rm mod}^{\hspace{.5pt}b\hspace{-2.5pt}}\La$. Thus, Statements (1) and (3) hold in this case. The proof of the theorem is completed.

\begin{Remark}\label{rlb-ass} In case $\La$ is a finite dimensional algebra, Theorem \ref{ass_fdm_cat}(3) is known to some specialists.

\end{Remark}

By Lemma \ref{llb_lrb} and Theorem \ref{ass_fdm_cat}, we have immediately the following result. 
\vspace{-1pt}

\begin{Cor}\label{ass_fdm_cor}

Let $\La=kQ/R$ be a locally semiprimary algebra, where $Q$ is a locally finite quiver and $\mk I$ is a relation ideal in $kQ$. If $\La$ is locally $($left, right$\hspace{.5pt})$  bounded, then ${\rm mod}^{\hspace{.5pt}b\hspace{-2.5pt}}\La$ has almost split sequences $($on the left, on the right$\hspace{.5pt})$.

\end{Cor}

\begin{Remark} Let $\La=kQ/\mk I$ be finite dimensional, where $Q$ is a finite quiver and $\mk I$ is a relation ideal in $kQ$ containing $(kQ^+)^n$ for some $n\ge 2$. Let $\tilde{\hspace{-2.5pt}\La}=k \tilde{Q}/\tilde{\mk I}$ such that there exists a quiver-covering $\phi: \tilde{Q}\to Q;$ see \cite[(4.1)]{BauL} and $\tilde{\mk I}$ is the pre-image of $\mk I$ under the $k$-linearly induced map $\phi': k \tilde{Q}\to kQ$. Then, $\phi'$ induces a covering $\psi: \tilde{\hspace{-2.5pt}\La} \to \La$ as defined in \cite[(3.1)]{BoG}. It is easy to check that $\tilde{\hspace{-2.5pt}\La}$ is locally semiprimary and locally bounded. By Corollary \ref{ass_fdm_cor}, $\mmod^{\hspace{.5pt}b\hspace{.5pt}}\tilde{\hspace{-2.5pt}\La}$ has almost split sequences, and by Theorem \ref{ass_fdm_cat}, the almost split sequences in $\mmod^{\hspace{.5pt}b\hspace{.5pt}}\tilde{\hspace{-2.5pt}\La}$ are almost split in $\Mod\hspace{2.5pt}\tilde{\hspace{-2.5pt}\La}.$

\end{Remark}

\begin{Exam}\label{ex_ass_fdim_cat} Consider $\La=kQ/\mk I$, where $Q$ is the quiver \vspace{-3pt} 
$$\begin{tikzpicture}[-{Stealth[inset=0pt,length=4.5pt,angle'=35,round,bend]},scale=.5]
\node at (-8,-0.02){$\cdots$}; 
\draw (-7.5,0) -- (-6,0);
\node at (-5.6,-0.05) {$a_3$};
\draw (-5.2,0) -- (-3.7,0);
\node at (-3.3,-0.05){$a_2$};
\draw  (-2.85,0) -- (-1.35,0);
\node at (-0.95,-0.05){$a_1$};
\node at (-2.2,0.3){\footnotesize $\alpha$};
\draw  (-0.55,-0.15) -- (.95,-0.8);
\node at (0,-0.9) {\small $\gamma$};
\draw  (-0.55, 0.15) -- (.95,0.8);
\node at (0.05, 0.8) {\footnotesize $\beta$};
\node at (1.35, 0.8) {$c_1$};
\node at (1.35,-0.9) {$c_2$};
\draw  (1.75, -0.9) -- (3.25,-0.07);
\draw  (1.75, 0.8) -- (3.25, 0.07);
\node at (2.65, -0.85) {\footnotesize  $\zeta$};
\node at (2.65, 0.8) {\small $\delta$};
\node at (3.6,-0.05) {$b_1$};
\draw (3.9,0) -- (5.4,0);
\node at (4.5,0.4){\footnotesize $\eta$};
\node at (5.75,-0.05) {$b_2$};
\draw (6.15,0) -- (7.65,0);
\node at (8,-0.05) {$b_3$};
 \draw (8.45,0) -- (9.95,0);
\node at (10.5,-0.02){$\cdots$};
\end{tikzpicture}\vspace{-8pt} $$ and $\mk{I}=\langle \beta \alpha, \gamma\alpha, \eta \delta, \zeta \eta \hspace{.8pt} \rangle.$ Clearly, $\La$ is neither locally left bounded nor locally right bounded. As seen in Example \ref{ex_ass_fp_cat}(1), the infinite dimensional indecomposable modules in ${\rm mod}^{-\hspace{-3pt}}\La$ are the injective modules $I_{a_i}$ with $i\ge 1$, \vspace{-.5pt} and the infinite dimensional indecomposable modules in ${\rm mod}^{+\hspace{-3pt}}\La$ are the projective modules $P_{b_j}$ \vspace{-1.5pt} with $j\ge 1$. By Theorem \ref{ass_fdm_cat}, ${\rm mod}^{\hspace{.5pt}b\hspace{-2.5pt}}\La$ has almost split sequences. Observe that ${\rm mod}^{\hspace{.5pt}b\hspace{-2.5pt}}\La$ has neither enough projective objects nor enough injective objects.

\end{Exam}


\vspace{5pt}

\begin{thebibliography}{99}

\vspace{1pt}


\bibitem{AnF} {\sc F. W. Anderson and K. R. Fuller}, ``Rings and Categories of Modules", Grad. Texts in Math. 13 (Springer-Verlag, New York, 1973).


\bibitem{Aus1} {\sc M. Auslander}, ``Coherent functors,'' Proc. Conf. Categorical Algebra (La Jolla, Calif., 1965) (Springer-Verlag, New York, 1966) 189-231.

\bibitem{A1} {\sc M. Auslander}, ``Representation theory of artin algebras I,'' 
Comm. Algebra 1 (1974) 177-268.

\bibitem{A2} {\sc M. Auslander}, ``Functors and morphisms determined by objects,'' Lect. Notes in Pure Appl. Math. 37 (Marcel Dekker, New York, 1978) 1-244.

\bibitem{Aus2} {\sc M. Auslander}, ``A survey of existence theorems for almost split sequences,'' London Math. Soc. Lecture Notes Ser. 116 (Cambridge University Press, Cambridge, 1986) 81-89.

\bibitem{Aus3} {\sc M. Auslander}, ``Almost split sequences and algebraic geometry,'' London Math. Soc. Lecture Notes Ser. 116 (Cambridge University Press, Cambridge, 1986) 165-179.

\bibitem{AIR} {\sc M. Auslander and I. Reiten}, ``Stable equivalence of dualizing R-varieties", Adv. Math. 12 (1974) 306-366. 

\bibitem{AuR2} {\sc M. Auslander and I. Reiten}, ``Representation theory of artin algebras III," Comm. Algebra 3 (1975) 239-294.

\bibitem{MAIR} {\sc M. Auslander and I. Reiten,} ``Representation theory of artin algebras IV," Comm. Algebra 5 (1977) 443-518.

\bibitem{ARS} {\sc M. Auslander, I. Reiten and Smal\o,} ``Representation Theory of Artin Algebras," Cambridge Stud. Adv. Math. 36 (Cambridge University Press, Cambridge, 1997).


\bibitem{BauL} {\sc R. Bautista and S. Liu}, ``Covering theory for linear categories with application to derived categories,'' J. Algebra 406 (2014) 173-225. 

\bibitem{BaL} {\sc R. Bautista and S. Liu}, ``The bounded derived categories of an algebra with radical squared zero,'' J. Algebra 482 (2017) 303-345. 

\bibitem{BLP} {\sc R. Bautista, S. Liu and C. Paquette}, ``Representation theory of strongly locally finite quivers," Proc. London Math. Soc. 106 (2013) 97-162.


\bibitem{BoG} {\sc K. Bongartz and P. Gabriel},``Covering spaces in representation theory", Invent. Math. 65 (1982) 331-378.


\bibitem{CKQ} {\sc W. Chin, M. Kleiner and D. Quinn,} ``Local theory of almost split sequences for comodules," Ann. Univ. Ferrara Sez. VII Sci. Mat. 51 (2005) 183-196.

\bibitem{CDi} {\sc W. Chin and D. Simson,} ``Coxeter transformation and inverses of Cartan matrices for coalgebras," J. Algebra 324 (2010) 2223-2248.

\bibitem{CEi} {\sc H. Cartan and S. Eilenberg}, ``Homological Algebra," (Princeton University Press, Princeton, New Jersey, 1956).

\bibitem{DSk} {\sc P. Dowbor and A. Skowro\'{n}ski,} ``Galois coverings of representation-infinite algebras," Comment. Math. Helv. 62 (1987) 311-337.

\bibitem{GaR}{\sc P. Gabriel and A. V. Roiter,} ``Representations of finite dimensional algebras,'' Encyclopedia Math. Sci. 73 (Springer-Verlag, Berlin, 1992).



\bibitem{Gup} {\sc E. Gupta}, ``Auslander-Reiten theory and string modules under the locally finite-dimensional setting,'' Master's thesis (Universit\'e de Sherbrooke, 2023).

\bibitem{Ha2} {\sc D. Happel}, ``Triangulated categories in the representation theory of finite dimensional algebras,'' London Math. Soc. Lecture Note Ser. 119 (Cambridge University Press, Cambridge, 1988).


\bibitem{HJor} {\sc T. Holm and P. J\o rgensen,} ``On a cluster category of infinite Dynkin type, and the relation to triangulations of the infinity-gon," Math. Z. 270  (2012) 277-295.

\bibitem{OIHY} {\sc O. Iyama, H. Nakaoka and Y. Palu,} ``Auslander-reiten theory in extriangulated categories,'' Trans. Amer. Math. Soc., Series B, 11 (2024) 248-305.

\bibitem{Jo1} {\sc P. J\o rgensen,} ``Auslander-Reiten theory over topological spaces,'' Comment. Math. Helv. 79 (2004) 160-182.

\bibitem{Jor} {\sc P. J\o rgensen,} ``Auslander-Reiten sequences on schemes," 
Ark. Mat. 44 (2006) 97-103.

\bibitem{HKra} {\sc H. Krause,} ``A short proof for Auslander's defect formula,'' Linear Algebra Appl. 365 (2003) 267-270.

\bibitem{Krau} {\sc H. Krause,} ``Krull-Remak-Schmidt categories and projective covers,'' Expo. Math. 33 (2015) 535-549.

\bibitem{HKra1} {\sc H. Krause,} ``Auslander-Reiten duality for Grothendieck abelian categories,'' Trans. Amer. Math. Soc. 371 (2019) 2455-2472.

\bibitem{HKS} {\sc H. Krause and M. Saor\'in,} ``On minimal approximations of modules,'' Contemp. Math. 229 (Amer. Math. Soc., Providence, RI, 1998) 227-236.



\bibitem{Lam} {\sc T. Y. Lam}, ``A First Course in Noncommutative Rings,'' Grad. Texts in Math. 131 (Springer-Verlag, New York, 1991).

\bibitem{LeZ} {\sc H. Lenzing and R. Zuazua}, ``Auslander-Reiten duality for abelian categories,'' Bol. Soc. Mat. Mexicana 10 (2004) 169-177.

\bibitem{LW} {\sc G. Leuschke and R. Wiegand}, ``Cohen-Macaulay representations," Math. Surveys Monogr. 181 (American Mathematical Society, Providence, RI, 2012). 

\bibitem{LL25} {\sc Z. Lin and S. Liu,} ``Representation theory of graded algebras given by locally finite quivers,'' J. Algebra 683 (2025) 672-718.



\bibitem{LiN} {\sc S. Liu and H. Niu,} ``Almost split sequences in tri-exact categories,'' J. Pure Appl. Algebra 226 (2022) 1-31.

\bibitem{SLPC} {\sc S. Liu, P. Ng and C. Paquette,} ``Almost split sequence and approximation,'' Algebr. Represent. Ther. 16 (2013) 1809-1827.

\bibitem{LiP} {\sc S. Liu and C. Paquette}, ``Cluster categories of type A double infinity and triangulations of the infinite strip," Math. Z. 286 (2017) 197-222.


\bibitem{RMV} {\sc R. Martinez-Villa,} ``Graded self-injective and Koszul algebras," J. Algebra 215 (1999) 34-72.



\bibitem{IRMVB}{\sc I. Reiten and M. Van den Bergh,} ``Noetherian hereditary abelian categories satisfying Serre duality,'' J. Amer. Math. Soc. 15 (2002) 295 - 366. 

\bibitem{RS}{\sc K. W. Roggenkamp and J. W. Schmidt,} ``Almost split sequences for integral group rings and orders,'' Comm. Algebra 4 (1976) 893-917. 




\bibitem{VHW} {\sc H. J. Von H\"{o}hne and J. Waschb\"{u}sch}, ``Die struktur
n-reihiger algebren,'' Comm. Algebra 12 (1984) 1187-1206.

\bibitem{Wat}{\sc C. E. Watts,} ``A homology theory for small categories", Proc. Conf. Categorical Algebra (La Jolla, Calif., 1965) (Springer-Verlag, New York, 1966) 331-335.

\vspace{5pt}

\end{thebibliography}
\end{document}